\title{On the existence of a \\ Kazantzis--Kravaris / Luenberger 
observer}
\author{Vincent Andrieu%
\thanks{V. Andrieu
(Vincent.Andrieu@onera.fr) is with DPRS, ONERA, BP 72, 29 Ave de la Division 
Leclerc, 92322 Chatillon, France
        }%
        \and Laurent Praly%
\thanks{L. Praly (Laurent.Praly@ensmp.fr) is with the Centre Automatique et Syst\`{e}mes,
           \'Ecole des Mines de Paris, 35 Rue Saint Honor\'{e},
           77305 Fontainebleau, France
        }
}

\documentclass[final]{siamltex}

\usepackage{amsfonts,latexsym}

\def\RR{{\mathbb R}}
\def\CC{{\mathbb C}}

\def\hx{\widehat{x}}

\def\NN{{\mathbb N}}
\def\OR{\mathcal{O}}
\def\BR{\mathcal{B}}

\def\DR{\mathcal{D}}
\def\bX{\breve{X}}

\def\bsigma{{\breve{\sigma }\,}}

\def\bf{\breve{f}}
\def\Vf{V_{\mathfrak{f}}}
\def\Vb{V_{\mathfrak{b}}}
\def\Wb{W_{\mathfrak{b}}}
\def\gf{\gamma _{\mathfrak{f}}}
\def\gb{\gamma _{\mathfrak{b}}}

\def\cl{\mathtt{cl}}
\def\Re{{\mathtt{Re}}}

\def\Box{\vbox{\hrule height0.6pt\hbox{%
   \vrule height1.3ex width0.6pt\hskip0.8ex
   \vrule width0.6pt}\hrule height0.6pt
  }}

\def\xrond{{\mathchoice
{{\mbox{$\scriptstyle \mathcal{X}$}}}
{{\mbox{$\scriptstyle \mathcal{X}$}}}
{{\mbox{$\scriptscriptstyle \mathcal{X}$}}}
{{\mbox{$\scriptscriptstyle \mathcal{X}$}}}
}}

\catcode`\@=11
\def\downparenfill{$\m@th\braceld\leaders\vrule\hfill\bracerd$}
\def\overparen#1{\mathop{\vbox{\ialign{##\crcr\crcr
\noalign{\kern0.4ex}
\downparenfill\crcr\noalign{\kern0.4ex\nointerlineskip}
$\hfil\displaystyle{#1}\hfil$\crcr}}}\limits}
\catcode`\@=12

\begin{document}
\maketitle

\begin{abstract}
We state sufficient conditions for the existence, on a given open 
set, of the
extension, to non linear systems, of the Luenberger observer as it has
been proposed by Kazantzis and Kravaris. 
We prove it is sufficient to choose the dimension of the 
system, giving the observer, less than or equal to 2 + twice the 
dimension
of the state to be observed. We show that it is sufficient to
know only an approximation of the solution of a PDE, needed for the
implementation. We establish a link with high gain observers.
Finally we extend our
results to systems satisfying an unboundedness observability property.
\end{abstract}

\pagestyle{myheadings}
\thispagestyle{plain}
\markboth{V. ANDRIEU AND L. PRALY}%
{ON THE KAZANTZIS--KRAVARIS / LUENBERGER OBSERVER}

\section{Introduction}
We consider the system~:
\begin{equation}
\label{3}
\dot x\;=\; f(x)\quad ,\qquad y\;=\; h(x)
\end{equation}
with state $x$ in $\RR^n$ and output $y$ in $\RR^p$
and where the functions $f$ and $h$ are sufficiently smooth.
We are concerned with the problem of existence of an observer for $x$ 
from the 
measurement $y$.

In  a seminal paper \cite{Kazantzis-Kravaris}, Kazantzis and Kravaris
have proposed to extend to the 
nonlinear case the primary observer introduced by Luenberger in 
\cite{Luenberger}
for linear systems. Following this suggestion,
the estimate $\hx$ of
$x$ is obtained as the output of the dynamical system~:
\begin{equation}
\label{2}
\dot z \;=\;  A\,  z\;+\; B(y)
\quad ,\qquad 
\hx \;=\; T^{*}(z)\  ,
\end{equation}
with state $z$ (a complex matrix) in $\CC^{m\times p}$ and where $A$ is a
Hurwitz complex matrix
and $B$ and $T^{*}$ are sufficiently smooth functions.

In the following we state sufficient conditions on $f$ and $h$
such that we can find $(A,B,m)$ for which there exists
$T^{*}$ guaranteeing the convergence of $\hx$ 
to $x$.

To ease readability, we have divided the paper into two parts. In a 
first part, we introduce and state our main results which are proved 
in the second part. Our first result gives a sufficient condition 
on 
$f$, $h$, $A$ and $B$ implying the existence of $T^{*}$ providing 
an appropriate observer. This condition involves a partial 
differential equation whose solution should be injective. In our second result, we propose a set of 
assumptions guaranteeing the existence of a solution for this 
equation.
Our third and fourth results give 
two sufficient conditions implying the injectivity property of this 
solution. Our 
fifth result shows that an observer can already be 
obtained if we know only an appropriate approximation of this solution.
This latter result allows us  to propose a new insight in the standard high
gain observer. Finally we claim that all these statements can be
extended to the case where the system satisfies an unboundedness
observability property.
\par\vspace{1em}\noindent
\textit{Some notations~:}
We assume the functions $f$ and $h$ in (\ref{3}) are at least locally 
Lipschitz.
So, for each $x$ in $\RR^n$, there exists a unique solution  $X(x,t)$ to 
(\ref{3}), with $x$ as initial condition.

Given an open set $\OR $ 
of $\RR^n$, for each $x$ in $\OR $,
we denote by 
$(\sigma_{\OR }^-(x),\sigma_{\OR }^+(x))$ 
the maximal interval of definition of the solution $X(x,t)$
conditioned to take values in $\OR $.

{}For a set $S$, we denote 
by $\cl(S)$ its closure and by $S +\delta $ the open set~:
\\[0.5em]
\null \hfill $\displaystyle 
S +\delta \;=\;  \left\{x\in \RR^n:\,  \exists\,   \xrond \in 
S :\,  |x-\xrond| < \delta \right\}
\;=\; 
\bigcup_{x\in S } \BR _\delta (x )\  ,
$\hfill \null \\[0.1em]
where 
$\BR_{\delta }(x)$ denotes the open ball with center $x$ and radius $\delta 
$.

By $L_fV$ we denote the Lie derivative of $V$ when it makes sense, i.e.~:
$$
L_fV(x)\;=\; \lim_{h\rightarrow 0}\frac{V(X(x,h))-V(x)}{h}
\ .
$$
\indent
Finally, $B_{1m}$ denotes the following vector in $\RR^m$
\begin{equation}
\label{82}
B_{1m}\;=\; \renewcommand{\arraystretch}{0.8}\left(\begin{array}{ccc}
1&\ldots&1
\end{array}\right)^T
\end{equation}
\section{Results and comments}
\subsection{Existence of a
Kazantzis-Kravaris / Luenberger observer}
\label{SecSufCond}
In \cite{Kazantzis-Kravaris}, $m$, the row dimension of $z$, is chosen
equal to $n$, the dimension of $x$, and $T^*$ is the inverse $T^{-1}$ of a function $T$,
solution of the following partial differential equation~:
\begin{equation}
\label{1}
\frac{\partial T}{\partial x}(x) \,  f(x) \; = \; A\,  T(x) \; + \; 
B(h(x))
\  .
\end{equation}
The rationale for this equation, as more emphasized in
\cite{Krener-Xiao1} (see also \cite{Rapaport-Maloum}), is that, if $T$ is a
diffeomorphism satisfying (\ref{1}), then the change of coordinates~:
\begin{equation}
\label{4}
\zeta \;=\; T(x)
\end{equation}
allows us to rewrite the dynamics (\ref{3}) equivalently as~:
$$
\dot \zeta  \;=\;  A\,  \zeta \;+\; B(h(T^{-1}(\zeta )))
\quad ,\qquad 
y\;=\; h(T^{-1}(\zeta ))\  .
$$
We then have~:
$$
\dot{\overparen{z-\zeta }}\;=\; A\,  (z-\zeta )
\  .
$$
$A$ being Hurwitz, $z$
in (\ref{2}) is the state of an asymptotically convergent observer of
$\zeta = T(x)$. Then, if the function $T^*=T^{-1}$ is uniformly 
continuous,
$\hx =T^{*}(z)$ is  an asymptotically convergent observer of
$x  =T^{*}(\zeta )= T^{*}(T(x))$.

This way of finding the function $T^*$ has motivated active research
on the problem of the existence of an analytic
and invertible solution to (\ref{1})
(see \cite{Kazantzis-Kravaris,Krener-Xiao1}
for instance). But, it turns out that having a (weak) solution to (\ref{1}) which is
only continuous and uniformly injective is already sufficient.
By sufficiency we mean, here, that an observer is
appropriate if we have convergence to zero of the observation error
associated to any solution which remains in 
a given open set $\OR$. For 
the latter, we need~:
\begin{definition}[Completeness within $\OR$]~$\,  $
The system (\ref{3}) is forward (\mbox{\rm resp.} backward) complete
within $\OR$ if we have the implication, for each $x$ in $\OR$,
\begin{equation}
\label{83}
\sigma _{\OR}^+(x)\; <\; +\infty \qquad \Longrightarrow\qquad 
\sigma _{\OR}^+(x) \; <\; \sigma _{\RR^n}^+(x)
\  .
\end{equation}
\end{definition}
In other words, completeness within $\OR$ says
that any solution $X(x,t)$ which exits $\OR$ in finite time must cross
the boundary of $\OR$ (at a finite distance).
An usual case where this property holds is when $f$ has an at most linear
growth within $\OR $.
\begin{theorem}[Sufficient condition of existence of an 
observer]
\label{ExistenceObs}
Assume,
the system (\ref{3}) is forward complete 
within $\OR$ and
there exist an integer $m$, a Hurwitz complex
$m\times m $ matrix $A$ and functions
$T:\cl (\OR)\rightarrow\CC^{m\times p}$, continuous,
$B:\RR^{p}\rightarrow\CC^{m\times p}$, continuous,
and
$\rho$, of class $\mathcal{K}_\infty$, satisfying~:
\begin{eqnarray}
\label{EDF}
&\displaystyle 
L_fT(x) = AT(x) + B(h(x))
\qquad \forall x \in \OR\  ,
\\
\label{UnifInj}
&\displaystyle 
|x_1-x_2|  \leq \rho(|T(x_1)-T(x_2)|)
\qquad \forall (x_1,x_2) \in \cl (\OR)^2
\  .
\end{eqnarray}
Under these conditions,
there exists a continuous function $T^*:\CC^{m\times p}\rightarrow 
\cl (\OR)$ such that, for each $x$ in $\OR$ and $z$ in
$\CC^{m\times p}$ the (unique) solution $(X(x,t),Z(x,z,t))$ of~:
\begin{equation}
\label{SysObs}
\dot x \;=\;  f(x)
\quad ,\qquad
\dot z \;=\;  Az+B(h(x))
\end{equation}
is right maximally defined on $[0,\sigma^+_{\RR^n}(x))$. Moreover, 
we have the implication~:
\begin{equation}
\label{34}
\sigma_{\OR }^+(x)\;=\;  \sigma_{\RR^{n}}^+(x)
\qquad \Longrightarrow\qquad 
\lim_{t\rightarrow \sigma_{\RR^n}^+(x)}
\left|T^{*}(Z(x, z, t))-X(x,t)\right|\;=\; 0
\  . 
\end{equation} 
\end{theorem}
\begin{remark}
\begin{remunerate}
\item
With the forward completeness within $\OR$ (\ref{83}), the condition on the left in (\ref{34}) 
implies that the solution $X(x,t)$ never exits $\OR$ and so~:
\begin{equation}
\label{84}
\sigma_{\OR }^+(x)\;=\;  \sigma_{\RR^{n}}^+(x)\;=\; +\infty 
\  .
\end{equation}
\item
Theorem \ref{ExistenceObs} extends readily
to the case where a) $y$ is a scalar, b) the state $x$ can be decomposed in
$x=(\xi _1,\xi _2)$ and
satisfies~:
$$
{\dot \xi }_1\;=\; f_1(\xi _1,u)\;+\; h(\xi _2)
\quad ,\qquad 
{\dot \xi }_2\;=\; f_2(\xi _2)
\quad ,\qquad y\;=\; \xi _1
$$
and c) the function $B$ can be chosen linear. In this case the observer
is implemented as the reduced order observer~:
$$
\renewcommand{\arraystretch}{0.1}\begin{array}{c}
\dot{\overparen{z-By}}\;=\; Az \;+\;  Bf_1(y,u)
\quad ,\qquad 
\widehat{\xi }_2\;=\; T^*(z)
\  .
\\[0.2em]
\null 
\end{array}
$$
\end{remunerate}
\end{remark}
\par\vspace{1.5em}
Assuming we have a continuous function $T$  satisfying
(\ref{EDF}), to implement the observer, we have
to find a function $T^{*}$ satisfying~:
$$
|T^{*}(z) -x|\; \leq \; \rho ^*(|z-T(x)|)
\qquad
\forall (x,z) \in \OR\times \CC^{m\times p}
$$
for some function $\rho ^*$ of
class $\mathcal{K}_\infty$. As shown by Kreisselmeier and Engel in
\cite{Kreisselmeier-Engel}, such a function $T^{*}$ exists if $T$ is continuous 
and uniformly injective as prescribed by (\ref{UnifInj}).

In conclusion, a Kazantzis-Kravaris / Luenberger observer
exists mainly if we can find a continuous function $T$ solving 
(\ref{EDF}) and uniformly injective in the sense of (\ref{UnifInj}).
\subsection{Existence of $T$ solving (\ref{EDF})}
To exhibit conditions guaranteeing the existence of a function $T$ 
solution of (\ref{EDF}), we abandon the 
interpretation above of a change of coordinates (see (\ref{4}))
and come back to the
original idea in \cite{Luenberger} (see also \cite{Kazantzis-Kravaris}
and \cite{Astolfi-Praly}) of dynamic extension.
Namely, we consider the augmented system (\ref{SysObs}).
Because of its triangular structure and the fact that $A$ is Hurwitz, 
we may expect this system to have, at least maybe only locally, an 
exponentially attractive invariant manifold in the augmented $(x,z)$ 
space
which could even be described 
as the graph of a function as~:
$$\left\{(x,z)\in\RR^n\times\CC^{m\times p}\,  :\: z=T(x)\right\}.$$
In this case, the function $T$ would satisfy the following identity,
for all $t$ in the domain of 
definition of the solution $(X(x,t),Z((x,z),t))$ of (\ref{SysObs}) 
issued from $(x,z)$ (compare with \cite[Definition 5]{Rapaport-Maloum}),
$$
T(X(x,t))\;=\; Z((x,T(x)),t)
\  ,
$$
or equivalently~:
\begin{equation}
\label{61}
T(X (x,t))\;=\; \exp(At)\,  T(x)
\;+\; \int_0^t \exp(As) B(h(X (x,s))) ds
\  .
\end{equation}
From this identity, (\ref{EDF}) is obtained by derivation with respect to $t$.
But, since we need (\ref{EDF}) to hold only on $\OR$,
from (\ref{61}), it is sufficient that $T$ satisfies~:
$$
T(x)\;=\; \exp(-At)\,  T(\bX (x,t))
\;-\; \int_0^t \exp(-As) B(h(\bX (x,s))) ds
\  ,
$$
where $\bX (x,s)$ is a solution of the modified system~:
\begin{equation}
\label{ModifSyst}
\dot x \;=\;  \bf(x)\;=\; \chi(x)\,  f(x)
\end{equation}
where $\chi:\RR^n\rightarrow\RR$ is an arbitrary locally Lipschitz
function satisfying ~:
\begin{equation}
\label{Chi}
\chi(x)\;=\; 1 \quad   \mathrm{ if } \quad   x\in \OR 
\quad ,\qquad 
\chi(x)\;=\; 0 \quad   \mathrm{ if } \quad   x\notin  \OR + \delta _u\ ,
\end{equation}
for some positive real number $\delta _u$. So, as standard in
the literature on invariant manifolds,
by letting $t$ go to $-\infty $, we get the 
following candidate expression for $T$~:
\begin{equation}
\label{7}
T(x)\;=\; \int_{-\infty }^0\exp(-As) B(h(\bX (x,s))) ds
\  .
\end{equation}
The above non rigorous reasoning can be made correct as follows~:

\begin{theorem}[Existence of $T$]
\label{Existence}
Assume the existence of a 
strictly positive real number $\delta _u$ such that
the system (\ref{3}) is backward complete 
within $\OR + \delta _u$.
Then, for each Hurwitz complex
$m\times m $ matrix $A$,
we can find a $C^1$ function $B:\RR^p\rightarrow\CC^{m\times p}$
such that the function $T:\cl(\OR)\rightarrow\CC^{m\times p}$,
given by (\ref{7}), is continuous and satisfies (\ref{EDF}).
\end{theorem}
\begin{remark}
All what is needed here about the function $B$
is that it guarantees that the 
function $t\mapsto |\exp(-At)B(h(\bX(x,t)))|$ is exponentially 
decaying with $t$ going to $-\infty $.
 So in particular (see Remark \ref{rem2}) when 
$\cl(\OR)$ is bounded, $B$ can be chosen simply as a linear function.
\end{remark}
\par\vspace{1em}
Approaching the problem from another perspective,
Kre\-isselmeier and Engel have introduced in 
\cite{Kreisselmeier-Engel}
this same expression (\ref{7}) (but with $X$ instead of $\bX$ and $B$ 
the identity function). Another link 
between \cite{Kazantzis-Kravaris} and
\cite{Kreisselmeier-Engel} has been established in
\cite{Krener-Xiao3}.
\subsection{$T$ injective}
Assuming now we have at our disposal the continuous function $T$, we
need to make sure that it is injective, if not uniformly injective as 
specified by (\ref{UnifInj}).
Here is where observability enters the game. Following 
\cite{Luenberger},
in \cite{Kazantzis-Kravaris,Krener-Xiao1}, when $m=n$, observability 
of the first order approximation at an equilibrium together with an 
appropriate choice of $A$ and $B$ is shown to imply injectivity of 
the solution $T$ of (\ref{1}) in a neighborhood of this equilibrium.
In \cite{Kreisselmeier-Engel}, uniform injectivity of $T$ is obtained 
under the following two assumptions~:
\begin{remunerate}
\item
The past output path $t\mapsto h(X(x,t))$ is uniformly injective
in $x$ with the set of past output paths equipped with an
exponentially weighted $L^2$-norm.
\item
The system (\ref{3}) has finite complexity, i.e.~there
exists a (finite) number $M$ of piecewise continuous function
$\phi_i$ in $L^2(\RR_-;\RR^p)$ and a strictly positive real number
$\delta $ such that we have, for each pair $(x_1, x_2)$ in $\OR^2$,
\\$\null\quad \displaystyle 
\sum_{i=1}^M
\left[\int_{-\infty }^0 \hskip -9pt \exp(-\ell s)\phi_i(s)^T
[h(X(x_1,s))-h(X(x_2,s))] ds\right]^2
$\hfill \null \\\null\hfill$\displaystyle
\geq \; 
\delta \,  \int_{-\infty }^0 \hskip -9pt\exp(-2\ell 
s)|h(X(x_1,s))-h(X(x_2,s))|^2 ds
\  .\quad  $
\end{remunerate}
\par\vspace{1.5em}
Our next result states, that, with the only assumption that
the past output path $t\mapsto h(X(x,t))$ is injective in $x$, it is 
sufficient to choose $m=n+1$ generic complex
eigen values for $A$ to get $T$ injective. 
The specific injectivity condition we need is~:
 \begin{definition}[Backward $\OR $-distinguishability]
\label{as4}
There exists two strictly positive real number $\delta _\Upsilon < \delta _d$  such that, 
for each
pair of distinct points $x_1$ and $x_2$ in $\OR +\delta _\Upsilon$,
there exists a time $t$, in $(
\max\{\sigma_{\OR +\delta _d}^-(x_1),\sigma_{\OR +\delta _d}^-(x_2)\}
\,  ,\,  0]$,
such that we have~:
$$
h(X(x_1,t))\; \neq \; h(X(x_2,t)) 
\  .
$$
\end{definition}
\indent
This distinguishability assumption says that the present state $x$
can be distinguished from other states in $\OR +\delta _\Upsilon$ by 
looking at the past output path restricted to the negative time 
interval where the solution $X(x,t)$ is in $\OR +\delta _d$.
\begin{theorem}[Injectivity]
\label{theo1}
Assume the system (\ref{3}) is backward complete
within $\OR + \delta_u$
and
backward $\OR$-distinguishable
with the corresponding $\delta _d$ in $(0,\delta _u)$.
Assume also the existence of an injective $C^1$ function $b:\RR^p\rightarrow\CC^p$,
a continuous function $M:\OR +\delta _\Upsilon\rightarrow\RR^+$,
and a negative real number $\ell$ such that, for each $x$ in 
$\OR +\delta _\Upsilon$,
the two functions
$t\mapsto \exp(-\ell t)b(h(\bX(x,t)))$ and
$t\mapsto \exp(-\ell t)\frac{\partial b\circ h \circ \bX}{\partial x}(x,t)$ 
satisfy, for each $t$ in $\displaystyle (\bsigma^-_{\RR^n}(x),0]$,
\begin{equation}
\label{65}
|\exp(-\ell t)b(h(\bX(x,t)))|\;+\; 
|\exp(-\ell t)\frac{\partial b\circ h \circ \bX}{\partial x}(x,t)|
\; \leq \; M(x)
\  ,
\end{equation}
where again $\bX$ is a solution of (\ref{ModifSyst}),
but this time with the function $\chi$ satisfying~:
\begin{equation}
\label{Chi2}
\chi(x)\;=\; 1 \quad  \mathrm{ if } \quad  x\in \OR + \delta _d 
\quad ,\qquad 
\chi (x)\;=\; 0 \quad  \mathrm{ if } \quad x\notin  \OR + \delta _u\ . 
\end{equation}
Under these conditions, there exists a subset $S$ of $\CC^{n+1}$ of 
zero Lebesgue measure
such that the function $T: \cl(\OR)\rightarrow\CC^{(n+1)\times p}$
defined, with the notation (\ref{82}), by~:
\begin{equation}
\label{62}
T(x)  \;=\;
\int_{-\infty}^0 \exp(-As) B_{1m}\,  b(h(\bX(x,s)))ds
\  ,
\end{equation}
is injective provided $A$ is a diagonal matrix with
$n+1$ complex eigen values $\lambda _i$ arbitrarily chosen 
in $\CC^{n+1}\setminus S$ and with real part strictly 
smaller than $\ell$.
\end{theorem}
\begin{remark}
\begin{remunerate}
\item
Condition (\ref{65}) holds for instance if
$f$, $h$ and $b$ have bounded derivative on $\cl(\OR +\delta _\Upsilon)$
(see \cite{Krener-Xiao3}).
\item
Theorem \ref{theo1} gives injectivity, not uniform injectivity. As already mentioned, if 
$\OR$ is bounded, continuity and the former imply the latter.
\end{remunerate}
\end{remark}
\par\vspace{1em}
Following Theorem \ref{theo1}, for any generic choice of $n+1$ complex eigenvalues 
for the matrix $A$, the function $T$ given 
by (\ref{62}) (or equivalently (\ref{7})) is injective. This says 
that the (real) row dimension of $z$ is
$m=2n+2$. It is a well known fact in observer theory that it is 
generically sufficient to extract $m = 2n+1$ pieces of 
information from the output path (with $h$ generically chosen)
to observe a state of dimension $n$  (see for instance 
\cite{Aeyels,Takens,Gauthier-Hammouri-Kupka,Coron,Sontag}). 
It can be understood from the adage that, the relation
$T(x_1)=T(x_2)$
between the two states $x_1$ and $x_2$ in $\RR^n$, i.e.~for $2n$
unknowns, has generically the unique trivial solution $x_1=x_2$ if we 
have 
strictly more than $2n$ equations, i.e.~$T(x)$ has strictly
more than $2n$ components.
\subsection{Injectivity in the case of complete observability}
Another setup where injectivity can be obtained is 
when we have complete observability. Namely
we can find a row dimension $m$ and a function
$b:y\in\RR^p\mapsto b(y) = (b_1(y),\dots,b_p(y))\in \RR^p$
so that
the following function $H:\RR^n\rightarrow \RR^{m\times p}$
is injective when restricted to $\cl(\OR)$~:
\begin{equation}
\label{H}
H(x) = \left ( \begin{array}{lll}b_1(h(x))&\dots&b_p(h(x))\\
L_{f}b_1(h(x))&\ldots&L_{f}b_p(h(x))\\
\vdots&\vdots&\vdots\\
L_{f}^{m-1}b_1(h(x)) &\ldots & L_{f}^{m-1}b_p(h(x))
\end{array}\right )
\  .
\end{equation}
Here $L_f^ih$ denotes the $i$th iterate Lie derivative,
i.e. $L_f^{i+1}h=L_f(L_f^ih$).
Of course, for this to make sense, the functions $b$, $f$ and $h$
must be sufficiently smooth.
This setup has been popularized and studied 
in deep details by Gauthier and his coworkers (see 
\cite{Gauthier-Kupka}
and the references therein, see also \cite{Rapaport-Maloum}).
In particular, it is established in \cite{Gauthier-Hammouri-Kupka} 
that, when $p=1$, for any generic pair $(f,h)$ in (\ref{3}), it is
sufficient to pick $m = 2n+1$.

With a Taylor expansion of the output path at $t=0$, we see
that the injectivity of $H$ implies that the function
which associates the initial condition $x$ to the output path, 
restricted to a very small time interval, is injective. This property 
is 
nicely exploited by observers with very fast dynamics as high gain 
observers (see \cite{Gauthier-Kupka}).
Specifically, we have~:
\begin{theorem}[Injectivity in the case of complete observability]
\label{theo2}
Assume the existence of a sufficiently smooth function $b:\RR^p\rightarrow\CC^p$
such that, for the function $H$
defined in (\ref{H}),
there exist a positive real number $L$ and a class $\mathcal{K}_\infty$ function
$\rho$ such that we have~:
\begin{eqnarray}
\label{LipCond2}
\null \hskip 2em
|L_{f}^mb(h(x_1)-L_{f}^mb(h(x_2))|&\leq &L|H(x_1)-H(x_2)|
\quad  \forall (x_1,x_2) \in \cl(\OR)^2\  ,
\\
\label{UnifInjH}
|x_1-x_2|& \leq & \rho(|H(x_1)-H(x_2)|)
\qquad 
\forall (x_1,x_2) \in \cl(\OR)^2\  .
\end{eqnarray}
Then, for any diagonal Hurwitz complex
$m\times m $ matrix $A$,
with $m$ the row dimension of $H$, there exists a real number $k^*$ such
that, for any $k$ strictly larger than $k^*$, there exists a function
$T:\cl(\OR)\rightarrow\CC^{m\times p}$ which 
is continuous, uniformly injective and satisfies (see (\ref{82}))~:
\begin{equation}
\label{49}
L_fT(x)  \;=\;  kAT(x) + B_{1m}b(h(x))
\qquad \forall x \in \OR
\  .
\end{equation}
\end{theorem}
\subsection{Approximation}
Fortunately  for the implementation of the observer, knowing
a function $T$ satisfying (\ref{EDF}) only approximately
is sufficient. But, in this case, we have to modify the
observer dynamics.
\begin{theorem}[Approximation]
\label{ExistenceObsAp}
Assume the system (\ref{3}) is forward complete within $\OR$.
Assume also
the existence of an integer $m$,
a Hurwitz complex
$m\times m $ matrix $A$ and functions
$T_a:\cl (\OR)\rightarrow\CC^{m\times p}$, continuous,
$B:\RR^{p}\rightarrow\CC^{m\times p}$, continuous, and
$\rho$ of class $\mathcal{K}_\infty$,
such that~:
\begin{equation}
\label{UCAp}
|x_1-x_2| \;  \leq \;  \rho(|T_a(x_1)-T_a(x_2)|)
\qquad \forall (x_1,x_2) \in \cl (\OR)^2\  ,
\end{equation}
the function $L_{f}T_a$ is well defined on $\OR$
and the function $\mathfrak{E}:\cl(\OR)\rightarrow\CC^{m\times p}$
defined as~:
\begin{equation}
\label{Eq21}
\mathfrak{E}(x) = L_{f}T_a(x) -  [AT_a(x) +B(h(x))]
\qquad  \forall x \in \OR
\end{equation}
satisfies~:
\begin{equation}
\label{EDFAp}
|\mathfrak{E}(x_1) - \mathfrak{E}(x_2)|\;\leq\; N\, 
|T_a(x_1)-T_a(x_2)|
\qquad
\forall (x_1,x_2) \in \cl (\OR)^2\  ,
\end{equation}
where $N$ is a positive real number satisfying~:
\begin{equation}
\label{RelLP}
2N\,   \lambda _{\max}(P) \;  < \;  1
\  ,
\end{equation}
with $\lambda _{\max}(P)$ the largest eigenvalue of the Hermitian 
matrix
$P$ solution of~:
\begin{equation}
\label{57}
\overline{A}^\top P + PA = -I
\  .
\end{equation}
Under these conditions,
there exists a function $T^*_a:\CC^{m\times p}\rightarrow 
\cl (\OR)$ and a locally Lipschitz function 
$\mathfrak{F}:\CC^{m\times p}\rightarrow \CC^{m\times p}$
such that, for each $x$ in $\OR$ and $z$ in
$\CC^{m\times p}$ each solution $(X(x,t),Z(x,z,t)) $ of~:
\begin{equation}
\label{SysObsAp}
\dot x \;=\; f(x)
\quad ,\qquad 
\dot z \;=\;  Az +\mathfrak{F}(z) +B(h(x))
\end{equation}
is right maximally defined on $[0,\sigma^+_{\RR^n}(x))$. Moreover, 
we have the implication~:
\begin{equation}
\label{81}
\sigma_{\OR }^+(x)\;=\;  \sigma_{\RR^{n}}^+(x)
\qquad \Longrightarrow\qquad 
\lim_{t\rightarrow \sigma_{\RR^n}^+(x)}
\left|T_a^{*}(Z(x, z, t))-X(x,t)\right|\;=\; 0
\  . 
\end{equation}
\end{theorem}
\begin{remark}
\begin{remunerate}
\item
In (\ref{Eq21}), $\mathfrak{E}$ represents the error in (\ref{EDF}) given by the 
approximation $T_a$ of $T$. This error should not be too large in an 
incremental sense as specified by (\ref{EDFAp}) and (\ref{RelLP}).
This indicates 
that one way to approximate $T$ is to look for $T_a$ in a set of functions 
minimizing the $L^\infty $ norm on $\cl (\OR)$ of the gradient of the 
associated error $\mathfrak{E}$. In particular, in the case where $\OR$ is
bounded, it follows from 
Weierstrass Approximation Theorem that we
can always choose
a Hurwitz complex matrix $A$ and a linear function $B$ so that
the constraint (\ref{RelLP}) can be satisfied by restricting ourself
to choose the function $T _a$ as a polynomial in $x$.
\item
The function $\mathfrak{F}$ in the observer
   (\ref{SysObsAp}) can be chosen as any Lipschitz extension of
   $\mathfrak{E}\circ T_a^*$ outside $T_a(\cl (\OR))$.
  This is very similar to what is done in \cite{Rapaport-Maloum}
  where a constructive procedure for this extension is proposed.
  Fortunately, this Lipschitz extension is not needed in the case 
where the function
  $\mathfrak{E}$ satisfies~:
	$$
	|\mathfrak{E}(x_1) - \mathfrak{E}(x_2)|\leq
 \frac{N}{4}\, \rho^{-1}(|x_1-x_2|)\quad
	\forall (x_1,x_2) \in \cl (\OR)^2\  ,
$$
where $\rho $ is the function satisfying (\ref{UCAp}).
  In this case we take simply (see
		(\ref{LipCondAp}))~:
  $$
  \mathfrak{F}(z) = \mathfrak{E}(T_a^*(z))\quad \forall z \in 
\CC^{m\times p}.
$$
\end{remunerate}
\end{remark}
\par\vspace{1em}
The combination of Theorems \ref{theo2} and \ref{ExistenceObsAp} gives us a new
insight in the classical
high gain observer of order $m$ as studied in 
\cite{Gauthier-Hammouri-Kupka}
or \cite{Rapaport-Maloum} for instance. Specifically, we have~:

\begin{corollary}[High gain Observer]
\label{cor1}
Assume
the system (\ref{3}) is forward complete within $\OR$
and
there exist a sufficiently smooth function 
$b:\RR^p\rightarrow\CC^p$,
a class $\mathcal{K}_\infty$ function $\rho$
and a positive real number $L$ such that
(\ref{LipCond2}) and (\ref{UnifInjH}) hold with $H$ and $m$ given by 
(\ref{H}).
Under these conditions,
for any diagonal Hurwitz complex
$m\times m $ matrix $A$,
there exists a real number $k^*$ such
that,
for any real number $k$ strictly larger than $k^*$, 
there exist a function $T^*_a:\CC^{m\times p}\rightarrow 
\cl (\OR)$, left inverse on $T_a(\cl (\OR))$ of the function
$T_{a}:\cl(\OR)\rightarrow\CC^{m\times p}$
defined as~:
\begin{equation}
\label{87}
T_a(x)\;=\; -\sum_{i=1}^{m} (kA)^{-i} B_{1m}  L_{f }^{i-1}b(h(x))
\  ,
\end{equation}
and a function $\mathfrak{F}:\RR^{m\times p} \rightarrow \RR^{m\times p}$
such that, for each $x$ in $\OR$ and $z$ in
$\RR^{m\times p}$ each solution $(X(x,t),Z(x,z,t)) $ of~:
\begin{equation}
\label{SysObsAp2}
\dot x \;=\; f(x)
\quad ,\qquad 
\dot z \;=\;  kAz \;+\; \mathfrak{F}(z) \;+\;  B_{1m}\,  b(h(x))
\end{equation}
is right maximally defined on $[0,\sigma^+_{\RR^n}(x))$. Moreover, 
we have the implication~:
\begin{equation}
\label{Concl2}
\null \qquad \null \sigma_{\OR }^+(x) \;=\; \sigma_{\RR^{n}}^+(x)
\qquad \Longrightarrow\qquad 
\lim_{t\rightarrow \sigma_{\RR^n}^+(x)}
\left|T_a^{*}(Z(x, z, t))-X(x,t)\right|\;=\; 0
\  . 
\end{equation}
\end{corollary}
\begin{remark}
When $\OR$ is bounded and $H$ is injective, uniform injectivity
(\ref{UnifInjH}) and forward completeness 
within $\OR$ hold necessarily. Thus, in this case, we recover
\cite[Lemma 1]{Rapaport-Maloum}.
\end{remark}
\subsection{Extension to boundedness observability}
\label{sec6}
Completeness is a severe restriction. Instead, it is proved in \cite{Astolfi-Praly}
that a necessary condition for the existence of an observer providing
the convergence to zero of
the observation error 
within the domain of definition of the solutions is the forward unboundedness observability property.
\begin{definition}[Unboundedness observable within $\OR$]
The system (\ref{3}) is forward (\mbox{\rm resp.} backward) unboundedness observable 
within $\OR$ if
there exists a proper and $C^1$ function
$\Vf : \RR^n \rightarrow \RR_+$ (\mbox{\rm resp.} $\Vb : \RR^n \rightarrow \RR_+$) and
a continuous function $\gf :\RR^p\rightarrow \RR_+$ (\mbox{\rm resp.}
$\gb :\RR^p\rightarrow \RR_+$) such 
that~:
\begin{eqnarray}
\label{ForUnb}
L_f \Vf (x)& \leq & \Vf (x)\;+\; \gf (h(x))
\qquad \forall x\,  \in \,  \OR
\  ,
\\\nonumber
(\; \mbox{\rm resp.}\  
L_f \Vb (x)& \geq &
-\Vb (x)\;-\; \gb (h(x))
\qquad \forall x\,  \in \,  \OR \  .\; )
\end{eqnarray}
\end{definition}
Fortunately, all our previous results still hold if completeness
is replaced by unboundedness observability but provided\footnote{%
The interested reader will find in \cite{Andrieu-Praly} the precise
statements of the corresponding results.
}~:
\begin{remunerate}
\item
the observer is modified in~:
$$
\dot z \;=\;   \gamma(y)\left[Az+B(y)\right]
\quad ,\qquad
\hx \;=\; T^{*}(z)\  ,
$$
where $\gamma $ is a $C^1$ function satisfying~:
$$
\gamma(h(x))\geq  1 + \gf (h(x))\qquad 
(\; \mbox{\rm resp. and}\  \gamma(h(x))\geq  1 + \gb (h(x))\; )
\qquad \forall x \in \cl (\OR)\  .
$$
As suggested in \cite{Astolfi-Praly}, the introduction of $\gamma $ takes care of 
possible
finite escape time. This has nothing in common
with the objective of error dynamics linearization as considered in
\cite{Respondek-Pogromsky-Nijmeijer}.
\item
In most occurrences, e.g.~(\ref{EDF}), (\ref{ModifSyst}),
(\ref{H}), (\ref{49}), (\ref{Eq21}), \ldots, $f$ is replaced by $f_\gamma 
$ defined as~:
$$
f_\gamma(x) \;=\;  \frac{f(x)}{\gamma(h(x))}
\  .
$$
\end{remunerate}

\section{Proofs}
\subsection{Proof of Theorem \ref{ExistenceObs}}
Because of the triangular structure of the system (\ref{SysObs}),
for each $x$ in $\OR$ and $z$ in 
$\CC^{m\times p}$, the component $Z(x,z,t)$ of the corresponding
solution of this system is defined as long as
$h(X(x,t))$ is defined. So this solution $(X(x,t),Z(x,z,t))$
is right maximally defined on the same interval
$[0,\sigma^+_{\RR^n}(x))$ as $X(x,t)$, solution of (\ref{3}).

Let us now restrict our attention to points $x$ in $\OR$ satisfying
the condition on the left in (\ref{34}). In this case, with the
forward completeness within $\OR$, we have (\ref{84}).
On the other hand, from (\ref{EDF}) and (\ref{SysObs}), we obtain, for each $x$ in $\OR$, $z$ in $\CC^{m\times p}$
and $t$ in $[0,\sigma^+_{\RR^n}(x))$,
\begin{equation}
\label{85}
T(X(x,t))-Z(x,z,t) \;=\;  
\exp\left(A t\right)
\,  (T(x)-z)
\  .
\end{equation}
As $A$ is a Hurwitz matrix, this and (\ref{84}) yield~:
$$
\lim_{t\rightarrow \sigma_{\RR^n}^+(x)}
\left|Z(x, z, t)-T(X(x,t))\right|\;=\; 0
\  . 
$$
From this, the implication (\ref{34}) follows readily if there exist a continuous function
$T^*:\CC^{m\times p}\rightarrow \cl (\OR)$ and a
class $\mathcal{K}_\infty$
function $\rho^*:\RR_+\rightarrow\RR_+$ satisfying~:
\begin{equation}
\label{UnifCont}
|T^*(z)-x|\; \leq \; \rho^*(z-T(x)) \qquad \forall z \in \CC^{m\times 
p}\; ,\  
\forall x \in \cl (\OR)
\  ,
\end{equation}
To find such functions, we first remark,
as in \cite{Kreisselmeier-Engel}, that 
(\ref{UnifInj}) and completeness of $\CC^{m\times p}$ and $\RR^n$
imply that $T(\cl (\OR))$ is a closed subset of $\CC^{m\times p}$.
It follows that, for each $z$ in $\CC^{m\times p}$, the infimum,
in $x$ in $\cl (\OR)$,
of $|T(x)-z|$ is achieved by at least one point, denoted $T_p^*(z)$
(in $\cl (\OR)$).
This defines a function
$T_p^*:\CC^{m\times p}\rightarrow \cl (\OR)$
satisfying~:
\begin{eqnarray}
\label{36}
&\displaystyle T(T_p^*(z))\;=\; z\qquad \forall z \in T(\cl (\OR))
\  ,
\\
\label{79}
&\displaystyle 
|T(T_p^*(z))-z|\; \leq \; |T(x)-z|\qquad \forall z \in \CC^{m\times 
p}\; ,\  \forall x \in \cl (\OR)\  .
\end{eqnarray}
With (\ref{UnifInj}), (\ref{36}) implies that the restriction $T_p^*$
to $T(\cl (\OR))$ is continuous. Also,  with the
triangle inequality, (\ref{79})
gives, for each $z$ in $\CC^{m\times p}$ and each $x$ in $\cl (\OR)$,
\begin{equation}
\label{35}
|x-T_p^*(z)| 
\; \leq\; 
\rho (|T(x)-z| + |z-T(T_p^*(z))|)
\; \leq\; \rho(2|T(x)-z|)
\  .
\end{equation}
Now we build the function $T^*$ by smoothing out $T_p^*$.
For each $z$ in $\CC^{m\times p}$, let
$$
\epsilon(z) \;=\; 
\frac{1}{2}\inf_{x\in \cl ( \OR)}|T(x)-z|
\  .
$$
$\left(\BR_{\epsilon(z)} (z )\right)_{z\in \CC^{m\times p} \setminus 
T(\cl (\OR))}$ 
is a covering of the open set $\CC^{m\times p} \setminus T(\cl (\OR))$
by open subsets.
From Lindel\"{o}f Theorem (see \cite[Lemma 4.1]{Boothby} for
instance),
there exists a sequence $\{z_i\}_{i\in\NN}$ such that 
$\left\{\BR_{\epsilon(z_i)}( z_i)\right\}_{i\in \NN}$ is a 
countable and locally finite covering by open subsets of  
$\CC^{m\times p} \setminus T(\cl (\OR))$.
{}For each $x$ in $\cl (\OR)$, each $z_i$ in $\{z_i\}_{i\in\NN}$
and each $z$ in $\BR_{\epsilon(z_i)}( z_i)$, we have~:
$$
|z_i-z|
\;  < \;  \epsilon(z_i)
\; \leq \; \frac{1}{2}\,   |T(x)-z_i|
\; \leq \; \frac{1}{2}\left[|T(x)-z| + |z-z_i|\right]
\; \leq \; |T(x)-z|
\  .
$$
With (\ref{35}), this gives~:
$$
|x-T^*_p(z_i)| 
\; \leq \;  \rho(2 |T(x)-z_i|)
\; \leq \;  \rho(2 (|T(x)-z| + |z-z_i |))
\; \leq \;  \rho(4| T(x)-z|)
\  .
$$
{}From \cite[Theorem IV.4.4]{Boothby}, we know that there exists a 
countable set of $C^\infty$ functions $\{\phi_i\}_{i\in 
\NN}:\CC^{m\times p} \setminus T(\cl (\OR))\rightarrow [0,1]$ 
satisfying, for each $z$ in $\CC^{m\times p} \setminus T(\cl (\OR))$,
$$
\sum_i \phi_i(z) \;=\;  1
\qquad ,\qquad \qquad 
   \phi_i(z) \;=\;  0\quad \forall z \notin \BR_{\epsilon(z_i)}( z_i) 
\  .
$$
We define the function $T^*:\CC^{m\times p}\rightarrow T(\cl (\OR))$
as~:
$$
\renewcommand{\arraystretch}{1.5}
\begin{array}{rcll}
T^*(z) &=& \displaystyle \sum_i \phi_i(z)T_p^*(z_i) &   
\mathrm{ if}\  z \in \CC^{m\times p} \setminus T(\cl (\OR))
\  ,
\\
 &=& T_p^*(z) &     \mathrm{if}\  z \in T(\cl (\OR))
\  .
\end{array}
$$

It is continuous when restricted to the open set $\CC^{m\times p} \setminus T(\cl
(\OR))$
and to the closed set $T(\cl (\OR))$.
Also, for each $z$ in $\CC^{m\times p}\setminus T(\cl (\OR) )$ and each
$x$ in $\cl (\OR)$, 
we get~:
\begin{eqnarray}
\nonumber
|T^*(z)-x|
\; = \; |\sum_{i} \phi_i(z) T_p^*(z_i)-x|
&\leq &\sum_{i} \phi_i(z)|T_p^*(z_i)-x|
\  ,\\\label{45}
&\leq &\sum_{i} \phi_i(z) \rho(4|z-T(x)|)
\ \leq \ \rho(4|z-T(x)|)
\  .\quad  \null 
\end{eqnarray}
And, for each $z$ in $T(\cl (\OR) )$ and each $x$ in $\cl (\OR)$, 
we get readily from (\ref{UnifInj}) and (\ref{36})~:
\begin{equation}
\label{46}
|T^*(z)-x|
\;=\; | T_p^*(z)-x|
\; \leq \; \rho (|T(T_p^*(z))-T(x)|)\;=\; \rho (|z-T(x)|)
\  .
\end{equation}
With (\ref{45}) and (\ref{46}), (\ref{UnifCont}) is established.
This proves also that $T^{*}$ is continuous on whole $\CC^{m\times p}$.
\subsubsection{Proof of Theorem \ref{Existence}}
With \cite[Corollary I.4.7]{Golubitsky-Guillemin}, we know
it exists a
locally Lipschitz function $\chi:\RR^n\rightarrow\RR$ satisfying 
(\ref{Chi}).
It follows that the function $\bf$ in (\ref{ModifSyst})
is locally Lipschitz. Thus, for each $x$ in $\RR^n$ there 
exists a unique solution 
$\bX(x,t)$ of (\ref{ModifSyst}), with initial condition $x$, 
maximally defined on 
$(\bsigma_{\RR^n }^-(x),\bsigma_{\RR^n }^+(x))$.
Moreover, backward completeness within $\OR + \delta _u$ of (\ref{3})
implies backward completeness of (\ref{ModifSyst}), i.e.~$\bsigma_{\RR^n}^-(x)=-\infty $.
Following \cite{Angeli-Sontag}, this implies the existence of
a proper and $C^1$ function
$\Vb : \RR^n \rightarrow \RR_+$ and
a continuous function 
$\gb :\RR^p\rightarrow \RR_+$ satisfying~:
\begin{equation}
\label{80}
L_{\bf} \Vb (x) \;   \geq \;  -\Vb (x) \;-\;  1
\qquad \forall x\,  \in \,  \RR^n\  .
\end{equation}
Let $\alpha$ be a strictly positive real number so that 
$A+\alpha  I$ is a Hurwitz matrix. We define the function
$\Wb :\RR^n\rightarrow\RR$ as~:
$$
\Wb (x) \;=\;  (\Vb (x)\;+\; 1)^\alpha 
\  ,
$$
With the help of Gronwall's Lemma, (\ref{80}) yields~:
\begin{equation}
\label{39}
\Wb (\bX(x,t))\leq \Wb (x) \exp(-\alpha t)
\qquad \forall x \in \RR^n\; ,\  \forall t \in 
(-\infty ,0] 
\  .
\end{equation}
Since $\Wb $ is a proper function and $h$ is continuous,
we can find a $C^1$ function
$\beta:\RR_+\rightarrow\RR_+$ of class $\mathcal{K}_\infty$
and a real number $\beta _0$
such that, for
each component $h_i$ of $h$, we have~:
$$
|h_i(x)|\; \leq \; \beta(\Wb (x))+\beta_0
\qquad \forall x \in \RR^n
\  .
$$
Let $\breve{\beta }:\RR_+\rightarrow \RR_+$ be the function defined as~:
$$
\breve{\beta }(w)\;=\; \sqrt{w}\;+\; \beta (w)\;+\; \beta _0
\  .
$$
This function is strictly increasing, $C^1$ on
$(0, +\infty)$ and its derivative $\breve{\beta }'$ satisfies~:
$$
\lim_{x\rightarrow 0}\breve{\beta }'(x) \;=\;  +\infty
\  .
$$
It admits an inverse $\breve{\beta }^{-1}:\RR_+\rightarrow \RR_+$
which satisfies~:
\begin{equation}
\label{38}
\breve{\beta }^{-1}(|h_i(x)|)\; \leq \; \Wb (x) \qquad 
\forall x \in \RR^n
\  .
\end{equation}
Moreover the function
$\eta \mapsto \frac{\eta \breve{\beta }^{-1}(|\eta |)}{|\eta |}$ is
$C^1$ on $\RR\setminus\{0\}$ and can be extended by continuity on $\RR$ as a 
$C^1$ 
injective function. So, with $p$ arbitrary vectors $b_j$ in $\RR^m$,
we define the function
$B:\RR^{p}\rightarrow\RR^{m\times p}$ as~:
$$
B(h) = \left (
\begin{array}{lll}
\frac{h_1 \breve{\beta }^{-1}(|h_1 |)}{|h_1 |} b_1& \dots &
\frac{h_p \breve{\beta }^{-1}(|h_p |)}{|h_p |} b_p
\end{array}\right )
\  .
$$
Since $A+\alpha I$ is a Hurwitz matrix, (\ref{39}), (\ref{38}) and the
backward completeness imply~:
\begin{remunerate}
\item
The existence of strictly positive real numbers $c_0$, $c_1$ and
$\varepsilon $ such that we have~:
\begin{eqnarray}
\null \hskip 2em|\exp(-As)B(h(\bX(x,s)))| &\leq &
c_0\,   |\exp(-As)|\,  \Wb (\bX(x,s)) 
\  ,
\\
\label{75}
& \leq & 
c_1\,  \Wb (x)\,  \exp(\varepsilon s)\qquad 
\forall (s,x)\in \RR_-\times \RR^n
\  .
\end{eqnarray}
\item
For each fixed $s$ in $\RR_-$, the function
$x\mapsto \exp(-As)B(h(\bX(x,s)))$ is continuous.
\end{remunerate}
\par\vspace{1em}\noindent
So Lebesgue dominated convergence Theorem (see
{\cite[Th\'{e}or\`{e}me (3.149)]{Deheuvels}} for instance)
implies that the following expression defines properly a continuous
function $T:\RR^n\rightarrow \CC^{m\times p}$~:
\begin{equation}
\label{T}
T(x)  \;=\;  \int_{-\infty}^0\exp(-As)\,  B(h(\bX(x,s)))ds
\  .
\end{equation}
Then,
for each $x$ in $\RR^n$ and for each $t$ in 
$(-\infty ,\bsigma_{\RR^n}^+(x))$, we get~:
\begin{eqnarray*}
T(\bX(x,t)) - T(x)&=& \int_{-\infty}^0\exp(-As)\,  
B(h(\bX(\bX(x,t),s)))ds
\;-\; T(x)
\  ,\\
&=& \exp(At)\,  \int_{-\infty}^{t}\exp(-Au)\,  B(h(\bX(x,u)))du\;-\; 
T(x)
\  ,\\
&=& (\exp(At)-I)\,   T(x) \,  +\,  
\exp(At)\,  \int_{0}^t\exp(-Au)B(h(\bX(x,u)))du
\  .
\end{eqnarray*}
Thus, we obtain, for all $x$  in $\RR^n$,
\begin{equation}
\label{52}
\chi(x)\,  L_fT(x) =
L_{\bf}T(x) =
\lim_{t\rightarrow 0}\frac{T(\bX(x,t)) - T(x)}{t} =
AT(x) + B(h(x))
\  .
\end{equation}
With (\ref{Chi}), this implies (\ref{EDF}) is satisfied.
\null \ \null\hfill\Box\par
\noindent
\begin{remark}
\label{rem2}
\begin{remunerate}
\item
For the case where $A$ is diagonalizable, with eigen value
$\lambda _i$,
and where the vectors $b_j$ are chosen so that the $p$ pairs $(A,b_j)$ are
controllable, our expression for $T$ gives for its $i$th component in 
the diagonalizing coordinates~:
\begin{equation}
\label{47}
T_i(x)\;=\; \int_{-\infty}^0\exp(-\lambda _i s)
B_i(h(\bX(x,s)))ds
\end{equation}
with~:
$$
B_i(h)\;=\; \left (
\begin{array}{lll}
\frac{h_1 \breve{\beta }^{-1}(|h_1 |)}{|h_1 |} b_{i1}& \dots &
\frac{h_p \breve{\beta }^{-1}(|h_p |)}{|h_p |} b_{ip}
\end{array}\right )
$$
where each $b_{ij}$ is non zero. Note that each of the $m$ rows of the function $B$ is an 
injective function from $\RR^p$ to $\RR^p$.
\item
If $\OR$ is bounded,  the function
$s\in \RR_-\mapsto h(\bX(x,s))\in\RR^p$ is a bounded function,
uniformly in $x$ in $\cl(\OR)$. It follows that the inequality
(\ref{75}) holds by choosing the function $\breve{\beta }^{-1}$ simply
as the identity function. This says that, in this case, the function
$B$ is linear.
\end{remunerate}
\end{remark}
\subsection{Proof of Theorem \ref{theo1}}
We first remark that backward $\OR $-dis\-tin\-gui\-sha\-bility property of the original
system (\ref{3}) implies the same property for the modified system (\ref{ModifSyst}).
Then we need the following Lemma.
\begin{lemma}[Coron]
\label{Coron}
Let $\Omega$ and $\Upsilon$ be open subsets of $\CC$ and $\RR^{2n}$ 
respectively.
Let $g: \Upsilon\times\Omega\rightarrow \CC^p$ be a
function which is holomorphic in $\lambda$ for each $x$ in $\Upsilon$
and $C^1$ in $x$ for each $\lambda $ in $\Omega $. 
If, for each pair $(x,\lambda )$ in $\Upsilon \times \Omega$ for 
which $g(x,\lambda )$ is zero 
we can find, for at least one of the $p$ components $g_j$ of $g$,
an integer $k$ satisfying~:
\begin{equation}
\label{Condg}
\null \qquad \null 
\frac{\partial^i g_j}{\partial \lambda^i}(x, \lambda) \;=\; 0
\qquad \forall i \in\{0,\ldots,k-1\}
\qquad ,\qquad \qquad 
\frac{\partial^{k} g_j}{\partial \lambda^{k}}(x, \lambda)\; \neq\;  0
\end{equation}
then the following set has zero Lebesgue measure in $\CC^{n+1}$~:
\begin{equation}
\label{58}
\null \qquad \null 
S\;=\; \bigcup_{x\in \Upsilon}
\left\{
 (\lambda_1,\dots,\lambda_{n+1})\in \Omega^{n+1} :\  
 g(x,\lambda _1) = \ldots=g(x,\lambda _{n+1})=0
\right\}\  .
\end{equation}
\end{lemma}
This result has been established by Coron in \cite[Lemma 3.2]{Coron} in a 
stronger form except for 
the very minor point that, here, $g$ is not $C^\infty $ in both $x$ and 
$\lambda $.
To make sure that this difference has no bad consequence and for the sake
of completeness, we give an ad hoc proof in appendix.

To complete the  proof of Theorem \ref{theo1}
all we have to do is to generate an appropriate function
$g$ satisfying all the required assumptions of this Lemma \ref{Coron}.

Let $\Omega $ and $\Upsilon$ be the following open subsets of $\CC$
and $\RR^{2n}$, respectively~:
$$
\Omega \;=\; \left\{\lambda \in \CC:\,  \Re(\lambda ) < \ell \right\}
\quad ,\qquad  
\Upsilon\;=\; 
\left\{x=(x_1,x_2)\in (\OR+\delta _\Upsilon)^2:\,  x_1\neq x_2\right\}
\  .
$$
By following the same arguments as in the proof
of Theorem \ref{Existence}, 
the backward completeness allows us to conclude~:
$$
\bsigma_{\RR^n}^-(x)=-\infty
\qquad \forall x \in \OR+\delta _\Upsilon
\  .
$$
Then, with (\ref{65}), we get,
for each $(x,\lambda,t )$ in
$(\OR+\delta _\Upsilon)\times\Omega \times(-\infty,0]$,
\begin{eqnarray*}
|\exp(-\lambda t)\,  b(h(\bX(x,t)))|
& \leq &
\exp([\ell - \Re(\lambda)] t)|\exp(-\ell t)b(h(\bX(x,t)))|
\  ,
\\
& \leq &
\exp([\ell - \Re(\lambda)] t)\,  M(x)
\;  .
\end{eqnarray*}
So Lebesgue dominated convergence Theorem implies that, for each 
fixed $\lambda $ in $\Omega $, the expression
$$
T_\lambda (x)\;=\; \int_{-\infty}^0\exp(-\lambda s)\,  b(h(\bX(x,s)))ds
$$
defines properly a continuous function $T_\lambda :\OR+\delta _\Upsilon
\rightarrow \CC^{p}$.
With similar arguments (see
{\cite[Th\'{e}or\`{e}me (3.150)]{Deheuvels}} for instance), 
with (\ref{65}), we can establish that  this function
is actually $C^1$.

Now, let $\DR T: (\OR+\delta _\Upsilon)^2\times \Omega \rightarrow \CC^p$ be the 
function defined as~:
\begin{equation}
\label{72}
\renewcommand{\arraystretch}{1.3}
\begin{array}[t]{rcl}
\DR T(x,\lambda)
&=& T_\lambda (x_1)- T_\lambda (x_2)
\  ,
\\
&=& \displaystyle \int_{-\infty}^0\exp(-\lambda s)[b(h(\bX(x_1,s)))-b(h(\bX(x_2,s)))]ds 
\  ,
\end{array}
\end{equation}
with $x=(x_1,x_2)$. It is 
$C^1$ in $x$ in $ (\OR+\delta _\Upsilon)^2$ for each $\lambda $ in $\Omega $.
Also, as proved in \cite[chap 19, p. 367]{Rudin}
with the help of Morera and Fubini Theorems, it is holomorphic in 
$\lambda $ in $\Omega $ for each $x$ in $ (\OR+\delta _\Upsilon)^2$. Moreover, 
since we have, for each $a$ in $(-\infty ,\ell)$ ,
$$
\int_{-\infty}^0 \exp(-2as)
|b(h(\bX(x_1,s)))-b(h(\bX(x_2,s)))|^2 \,  ds\; \leq \; 
\frac{M(x_1)^2+M(x_2)^2}{2(\ell -a)}\; <\; +\infty 
\  ,
$$
we can apply Plancherel Theorem to obtain,
for each $a$ in $(-\infty ,\ell)$ and each $x$ in $(\OR+\delta _\Upsilon)^2$,
\begin{equation}
\label{71}
\frac{1}{2\pi }\,  \int_{-\infty }^{+\infty }
|\DR T(x,a +is)|^2\,  ds
\;=\; 
\int_{-\infty}^0 \exp(-2a s)
|b(h(\bX(x_1,s)))-b(h(\bX(x_2,s)))|^2\,  ds
\  .
\end{equation}
Now, for $x$ in $\Upsilon$, with the distinguishability property, 
continuity with respect to time and injectivity of $b$ imply the existence of an open time
interval $(t_0,t_1)$ such that~:
$$
|b(h(\bX(x_1,s)))- b(h(\bX(x_2,s)))|\; > \; 0
\qquad \forall s \in (t_0,t_1)
\  .
$$
It follows with (\ref{71}) that we have~:
$$
\int_{-\infty }^{+\infty }
|\DR T(x,a +is)|^2\,  ds\; > \; 0
\  .
$$
This says that, for each $x$ in $\Upsilon$,
 the function $\lambda \mapsto \DR  T(x,\lambda)$ is not identically equal zero
on $\Omega $. Since it is holomorphic, this implies that,
 for each $(x,\lambda )$ in $\Upsilon\times \Omega $,
we can find, for at least one of the $p$ components $\DR T_j$ of
$\DR T$, an integer $k$ satisfying~:
$$
\frac{\partial^i \DR  T_j}{\partial \lambda^i}(x, \lambda) \;=\; 0
\qquad \forall i \in\{0,\ldots,k-1\}
\qquad ,\qquad \qquad 
\frac{\partial^k \DR  T_j}{\partial \lambda^k}(x, \lambda)\; \neq\;  0
\  .
$$

So we can invoke Lemma \ref{Coron} with $\DR$ as function $g$.
With (\ref{72}), it allows us to conclude that the following set $S$
has zero Lebesgue measure in $\CC^{n+1}$~:
$$
S\;=\; \bigcup_{(x_1,x_2)\in \Upsilon} \left\{
(\lambda_1,\dots,\lambda_{n+1})\in \Omega^{n+1} :\quad 
 T_{\lambda _i}(x_1) = T_{\lambda _i}(x_2)\quad  \forall i \in \{1,\ldots,n+1\}\right\}
\  .
$$
\subsection{Proof of Theorem \ref{theo2}}
Our first step consists in proposing a function $T$ solution of 
(\ref{49}). The definition (\ref{H}) of $H$ and the inequality 
(\ref{LipCond2}), give, for each pair $(x_1,x_2)$ in $\cl(\OR)^2$,
\begin{eqnarray*}
|L_{f}H(x_1)-L_{f}H(x_2)|  &\leq&
|H(x_1)-H(x_2)| + |L_{f}^mb(h(x_1))-L_{f}^mb(h(x_2))|
\  ,
\\
&\leq& (1+L)|H(x_1)-H(x_2)|
\  .
\end{eqnarray*}
Also (\ref{UnifInjH}) implies that, for each $Y$ in
$H(\cl (\OR))$, there exists a unique $x$ in $\cl (\OR)$ solution of 
$Y=H(x)$.
Hence we can define a Lipschitz function
$F:H(\cl(\OR)) \rightarrow \RR^{m\times p}$ satisfying~:
\begin{equation}
\label{53}
F(H(x)) \;=\; L_{f}H(x) \qquad \forall x \in \cl(\OR)\  .
\end{equation}
Furthermore, as in the proof of Theorem 
\ref{ExistenceObs},
continuity and uniform injectivity of the function $H$ on 
$\cl(\OR)$ as given by (\ref{UnifInjH})
imply that $H(\cl(\OR))$ is closed. Then it follows from
Kirszbraun's Lipschitz extension Theorem
(see \cite[Theorem 2.10.43]{Federer} for instance) that $F$
can be extended as a function $\breve{F}:\RR^{m\times p}\rightarrow\RR^{m\times p}$ 
satisfying~:
\begin{eqnarray}
\label{50}
|\breve{F}(Y_1)-\breve{F}(Y_2)|&\leq &(1+L)\,  |Y_1-Y_2|
\qquad \forall (Y_1, Y_2)  \in\RR^{m\times p} \times \RR^{m\times p}
\  ,
\\
\label{51}
\breve{F}(Y)&=&F(Y)\qquad \forall Y \in H(\cl (\OR))
\  .
\end{eqnarray}
Let $\mathfrak{Y}(Y,t)$ denote a solution of the following system on 
$\RR^{m\times p}$~:
$$
\dot Y = \breve{F}(Y)
\  .
$$
With (\ref{50}), such a solution is unique for each $Y$ in 
$\RR^{m\times p}$, defined on
$(-\infty ,+\infty )$ and satisfies, for some fixed matrix $Y_0$ in $\RR^{m\times p}$
and for each pair $(Y_1,Y_2 )$ in 
$\RR^{m\times p}\times \RR^{m\times p}$,
$$
\renewcommand{\arraystretch}{1}
\begin{array}{rcl}
|\mathfrak{Y}(Y_1,t)-Y_0|
&\leq &\displaystyle 
|Y_1-Y_0|
\;+\; \int_t^0
|\breve{F}(\mathfrak{Y}(Y_1,s))-\breve{F}(Y_0)|ds\;-\; |\breve{F}(Y_0)|\,  t
\  , 
\\
&\leq & \displaystyle 
|Y-Y_0|
\;+\; (1+L)\,  \int_t^0
|\mathfrak{Y}(Y,s)-Y_0|ds 
\;-\; |\breve{F}(Y_0)|\,  t
\  ,
\\
|\mathfrak{Y}(Y_1,t)-\mathfrak{Y}(Y_2,t)|&  \leq & \displaystyle 
 (1+L)\,  \int_t^0
|\mathfrak{Y}(Y_1,s)-\mathfrak{Y}(Y_2,t)|ds
\  .
\end{array}
$$
With Gronwall inequality, this gives, for all $t\leq 0$,
\begin{eqnarray}
\label{77}
|\mathfrak{Y}(Y,t)-Y_0|
& \leq &  \exp(-(1+L)t)
\left[|Y-Y_0|
+ \frac{|\breve{F}(Y_0)|}{1+L}\right]
\;-\;\frac{ \breve{F}(Y_0)}{1+L}
\  ,
\\\label{56}
\null \qquad \qquad \null |\mathfrak{Y}(Y_1,t)-\mathfrak{Y}(Y_2,t)| &\leq &  \exp(-(1+L)t)
\,  |Y_1-Y_2|
\  .
\end{eqnarray}
So, in particular, we have, for each $t\leq 0$ and  $Y$ in 
$\RR^{m\times p}$,
\begin{equation}
\label{78}
|\breve{F}(\mathfrak{Y}(Y,t))|\;\leq\;
(1+L)\left(\exp(-(1+L)t)
\left[|Y-Y_0|
+\frac{|\breve{F}(Y_0)|}{1+L}\right]\right)
\  .
\end{equation}
Hence, given any diagonal Hurwitz $m\times m$ matrix $A$, with eigen value $\lambda _i$,
for each real number $k\geq \frac{1+L}{-\max_i\{\Re(\lambda _i)\}} $,
we can properly define a continuous function
$\mathfrak{R}:\RR^{m\times p}\rightarrow \CC^m$ as~:
\begin{equation}
\label{67}
\mathfrak{R}(Y) \;=\;  \int_{-\infty }^0 \exp(-skA)\,  B_{1m}\,
\breve{F}( \mathfrak{Y}(Y,s))_m ds \  ,
\end{equation}
with the notation
(\ref{82}), and where $\breve{F}( Y)_m$ denotes the $m$th row of 
$\breve{F}( Y)$.
As for (\ref{52}), we can prove that we have~:
$$
L_{\breve{F}} \mathfrak{R}(Y)\;=\; kA\,  \mathfrak{R}(Y)\;+\; B_{1m}\,  
\breve{F}( Y)_m 
\qquad \forall Y\in\RR^{m\times m}
\  .
$$
But, with (\ref{H}), (\ref{53}) and (\ref{51}), this yields~:
\begin{equation}
\label{EqHG21}
L_F \mathfrak{R}(H(x))\;=\; kA\,  \mathfrak{R}(H(x))\;+\; B_{1m}\,  
L_{f}^m b(h(x))
\qquad \forall x\in \cl(\OR)
\  .
\end{equation}

Let now $T:\cl(\OR)\rightarrow\RR^n$ be the continuous function defined as~:
\begin{equation}
\label{55}
T(x) \;=\; (kA)^{-m}\mathfrak{R} (H(x)) - K^{-1}SH(x) \, ,
\end{equation}
with the notations~:
\begin{equation}
\label{88}
S \;=\;  \left ( 
\begin{array}{ccc}
\lambda_1^{-1} & \ldots & \lambda_1^{-m}\\
\vdots&\vdots&\vdots\\
\lambda_m^{-1}&\ldots & \lambda_m^{-m}
\end{array}
\right )
\quad ,\qquad K\;=\; \diag{(k,\ldots,k^m)}
\  .
\end{equation}
We want to show that $T$ is a solution of (\ref{49}).
We have~:
$$
K^{-1}SH(x) =  \left ( \begin{array}{c}
\sum_{i=1}^m (k\lambda_1)^{-i}L_{f}^{i-1}b(h(x))\\
\vdots\\
\sum_{i=1}^m (k\lambda_m)^{-i}L_{f}^{i-1}b(h(x))\\
\end{array}\right )
\  .
$$
Thus, for each $x$ in $\RR^n$, we get~:
\begin{equation}
\label{54}
\null \qquad K^{-1}SL_{f}H(x)
\;=\; 
kAK^{-1}SH(x) \;-\;  B_{1m}b(h(x)) \;+\; (kA)^{-m}B_{1m}L_{f}^mb(h(x))
\  .
\end{equation}
In view of (\ref{55}), it remains to compute the Lie derivatives of $(kA)^{-m}\mathfrak{R} 
(H(x))$.
From (\ref{Chi}), (\ref{H}), (\ref{53}) and (\ref{51}),
we get the identity~:
$$
\mathfrak{Y}(H(x),t)\;=\; H(\bX (x,t))
\qquad \forall t \in (\bsigma _\OR^-(x),\bsigma _\OR^+(x))
\qquad \forall  x \in \OR
\  .
$$
This gives readily, for all $t$ in
$(\bsigma _\OR^-(x),\bsigma _\OR^+(x))$ and $x$ in $\OR$,
$$
\mathfrak{R}(\mathfrak{Y}(H(x),t))\;-\; \mathfrak{R}(H(x))
\;=\; \mathfrak{R}(H(\bX (x,t))\;-\; \mathfrak{R}(H(x))
$$
and therefore~:
$$
L_F\mathfrak{R}(H(x))\;=\;
L_{f}\mathfrak{R}(H(x))
\qquad \forall x \in \OR
\  .
$$
By continuity this identity extends to $\cl(\OR)$.
So, with (\ref{EqHG21}), we get~:
$$
L_{f}\mathfrak{R}(H(x)) \;=\; kA\,  \mathfrak{R}(H(x))
\;+\; B_{1m}\,  L_{f}^m(b(h(x))) 
\qquad \forall x \in \cl(\OR).
$$
Consequently, with (\ref{55}) and (\ref{54}), we finally obtain, for each
$x$ in $\cl(\OR)$,
\begin{eqnarray*}
L_{f}T(x) &=&
(kA)^{-m}L_{f}\mathfrak{R} (H(x)) - K^{-1}SL_{f}H(x)
\  ,
\\
&=&
kA\left[(kA)^{-m}\,  \mathfrak{R}(H(x)) -K^{-1}SH(x)\right] + 
B_{1m}b(h(x))
\  ,
\\
&=&  kAT(x) + B_{1m}b(h(x))
\  .
\end{eqnarray*}
This proves that the function $T$ defined by (\ref{55}) is solution of 
(\ref{49}).

Our second step in this proof is to show that, by picking $k$ large enough,
the function $T$ 
given by (\ref{55}) is uniformly injective. 
To simplify the following notations, to a function $f$, we associate the function 
$\Delta f$ as follows~:
$$
\Delta f(x_1,x_2) \;=\;  f(x_1)\;-\; f(x_2)
\  .
$$
So, for instance, for each pair $(x_1,x_2)$ in $\OR^2$, we have~:
$$
T(x_1)-T(x_2) \;=\;  (kA)^{-m}\Delta (\mathfrak{R} \circ H)(x_1,x_2))
\;+\;  K^{-1}S\Delta H(x_1,x_2)
\  .
$$
With (\ref{50}), (\ref{67}) and (\ref{56}), we get, for each 
$(Y_1,Y_2)$ in
$\RR^{m\times p}\times \RR^{m\times p}$,
\begin{eqnarray*}
|\Delta \mathfrak{R} (Y_1,Y_2)| & \leq &
\int_{-\infty }^0\left| \exp(-skA)\,  B_{1m}\, 
\left[\breve{F}(\mathfrak{Y}(Y_1,s))_m-\breve{F}(\mathfrak{Y}(Y_2,s))_m\right]
\right|ds
\\
& \leq & (1+L)\,
\int_{-\infty }^0 |\exp(-skA)|\,  |B_{1m}|
\, \left[\mathfrak{Y}(Y_1,s)-\mathfrak{Y}(Y_2,s)\right|ds
\\
& \leq &(1+L)\,  |B_{1m}|\,  |Y_1-Y_2|\,  
\int_{-\infty }^0 \exp\left(-s\left[1+L+k\max_i\{\Re(\lambda 
_i)\}\right]\right) ds
\\
& \leq &
\frac{(1+L)|B_{1m}|}{-\left[1+L+k\max_i\{\Re(\lambda _i)\}\right]}
\,  |Y_1-Y_2 |
\  .
\end{eqnarray*}
This yields, for each pair $(x_1,x_2)$ in $\cl(\OR)^2$~:
\begin{eqnarray}
\nonumber
|T(x_1)-T(x_2)|
&\geq &
|K^{-1}S\Delta H(x_1,x_2)| \;-\;
|(kA)^{-m}\Delta (\mathfrak{R} \circ H)(x_1,x_2)| 
\  ,
\\\nonumber
&\geq &
\frac{|\Delta H(x_1,x_2)|}{|K|\,|S^{-1}|}
\;-\;
\frac{ |(kA)^{-m}|\,(1+L)|B_{1m}|}{-\left[1+L+k\max_i\{\Re(\lambda _i)\}\right]}
\,  |\Delta H(x_1,x_2)|
\  ,
\\\nonumber
&\geq &
k^{-m}\left(\frac{1}{|S^{-1}|}-
\frac{|A|^{-m}(1+L)|B_{1m}|}{-\left[1+L+k\max_i\{\Re(\lambda _i)\}\right]}
\right)
|H(x_1-H(x_2)|
\  .\quad  \null 
\end{eqnarray}
So, with (\ref{UnifInjH}), the function $T$ is uniformly injectivity on $\cl(\OR)$
for all $k$ large enough.
\subsection{Proof of Theorem \ref{ExistenceObsAp}}
Following the same arguments as in the proof of 
Theorem \ref{ExistenceObs}, 
continuity and uniform injectivity (\ref{UCAp}) of the function $T_a$ 
on $\cl(\OR)$ imply that $T_a(\cl(\OR))$ is closed and 
we can construct a continuous function $T_a^*:\CC^{m\times p}\rightarrow \cl (\OR)$ 
satisfying~:
\begin{equation}
\label{UnifCont2}
|T_a^*(z)-x|\; \leq \; \rho(4|z-T_a(x)|)
\qquad
\qquad 
\forall (x,z) \in \cl(\OR) \times \CC^{m\times p}
\  .
\end{equation}
This implies~:
\begin{equation}
\label{69}
T_a^{*}(T_a(x))\;=\; x
\qquad \forall x \in \cl(\OR) 
\  .
\end{equation}

Now, let us assume for the time being there exists a function 
$\mathfrak{F}:\CC^{m\times p}\rightarrow\CC^{m\times p}$ to be
used in (\ref{SysObsAp}) and satisfying~:
\begin{equation}
\label{LipCondAp}
|\mathfrak{E}(x)-\mathfrak{F}(z)| \;\leq \; N\, |T_a(x)-z| 
\qquad 
\forall (x,z) \in \cl(\OR) \times \CC^{m\times p}
\  .
\end{equation}
As a direct consequence, we get the inequality~:
$$
| \mathfrak{F}(z) | \; \leq \; N\,   |z| \;+\;  M 
\qquad z \in \CC^{m\times p}\  ,
$$
for some real number $M$ ($= |\mathfrak{E}(x_0)| + N|T_a(x_0)|$, 
with some arbitrarily fixed $x_0$  in $\cl(\OR)$). It follows that
the $z$ dynamics in the system (\ref{SysObsAp}) 
satisfy~:
$$
|\dot z| \;\leq\; (|A|+N)\,|z| \;+\; (M \,+\, |B(h(x)))|\, .
$$
Hence, for each $x$ in $\OR$ and $z$ in 
$\CC^{m\times p}$, the component $Z(x,z,t)$ of a solution
$(X(x,t),Z(x,z,t))$ of (\ref{SysObsAp}) is defined as long as
$h(X(x,t))$ is defined. So this solution 
is right maximally defined on the same interval 
$[0,\sigma^+_{\RR^n}(x))$ as $X(x,t)$ solution of (\ref{3}).

With (\ref{UnifCont2}), (\ref{81}) holds if we have~:
\begin{equation}
\label{Conv2}
\lim_{t\rightarrow \sigma^+_{\RR^n}(x)}T_a(X(x,t)) - Z(x,z,t) \;=\; 0
\qquad
\qquad 
\forall (x,z) \in \OR \times \CC^{m\times p}
\  .
\end{equation}
To establish this limit, we associate, to each pair $(x,z)$ in $\OR \times \CC^{m\times 
p}$, the matrix $e$ in $\CC^{m\times p}$~:
$$
e \;=\; T_a(x) - z
\  .
$$
With (\ref{Eq21}) and (\ref{SysObsAp}), we get~:
$$
\dot e \;=\;  Ae + \mathfrak{E}(x) - \mathfrak{F}(z)
\  .
$$
Let $U:\CC^{m\times p}\rightarrow \RR_+$ be the positive definite and 
proper function defined as~:
$$
U(e) \;=\; \sum_{i=1}^p \overline e_i^\top Pe_i\;,
$$
where $e_i$ denotes the $i^{th}$ column of $e$, $\overline{e_i}$ 
denotes its complex conjugate and $P$ is given by (\ref{57}).
Using (\ref{LipCondAp}) and completing the squares, we get~:
\begin{eqnarray}
\nonumber
\dot{\overparen{U(e)}} \;=\; 
\sum_{i=1}^p \left[-|e_i|^2 \;+\;
2\overline{e_i}^\top P
\left( \mathfrak{E}(x) - \mathfrak{F}((z)\right)_i\right]
&\leq &   -[1- 2N\lambda _{\max}(P)]\,  |e|^2
\  ,
\\\nonumber
&\leq &  - \frac{1- 2N\lambda _{\max}(P)}{\lambda _{\min}(P)}\,  U(e)
\  .
\end{eqnarray}
So, with (\ref{RelLP}), we have established the existence of a 
strictly positive real number $\varepsilon $ such that we have~:
\begin{equation}
\label{EqU1}
\dot{\overparen{U(e)}}\;\leq\; -  \varepsilon \,  U(e)
\qquad
\forall (x,z) \in \OR\times\CC^{m\times p}
\  .
\end{equation}
This implies, for all $t$ in $[0,\sigma _\OR^+(x))$ and $(x,z)$ in
$\OR\times\CC^{m\times p}$,
\begin{equation}
\label{86}
\null \qquad \quad \exp\left(-\varepsilon t\right)\,  U(e)\; \geq \;
U(E(x,z,t))\qquad \quad 
\left(\  \;=\;   U\left(T_a(X(x,t)) - Z(x,z,t)\right) \  \right) 
\  .
\end{equation}
With forward completeness within $\OR$ and the condition in the left 
of (\ref{81}) (see (\ref{84})), this implies (\ref{Conv2}) holds.

It remains to establish
the existence of 
a function $\mathfrak{F}:\CC^{m\times p}\rightarrow\CC^{m\times p}$ 
satisfying (\ref{LipCondAp}).
With (\ref{69}), we see that (\ref{EDFAp}) becomes~:
$$
|\mathfrak{E}(T_a^*(z_1))-\mathfrak{E}(T_a^*(z_2))| \;\leq \; N\,  |z_1-z_2| 
\qquad \forall (z_1,z_2) \in T_a(\cl(\OR))^2
\  .
$$
Thus, $\mathfrak{E}\circ T_a^*$ is a Lipschitz function
on the closed subset $T_a(\cl(\OR))$. From Kirszbraun's Lipschitz
extension Theorem, $\mathfrak{E}\circ T_a^*$ can be 
extended as a function
$\mathfrak{F}:\CC^{m\times p}\rightarrow\CC^{m\times p}$ 
satisfying~:
\begin{eqnarray*}
| \mathfrak{F}(z_1) \,-\, \mathfrak{F}(z_2) |
&  \leq  & N \, | z_1 \,-\, z_2|
\quad \forall (z_1,z_2) \in (\CC^{m\times p} )^2
\  ,
\\
\mathfrak{F}(z) & = & \mathfrak{E}(T^*(z)) \quad \forall z \in 
T_a(\cl(\OR))
\  .
\end{eqnarray*}
So, in particular, we get (\ref{LipCondAp}).
\subsection{Proof of Corollary \ref{cor1}}
Let $\lambda _i$ be the eigen values of a
given diagonal
Hurwitz complex
$m\times m $ matrix $A$.
With the notations (\ref{88}), the function $T_{a}:\cl(\OR)\rightarrow\RR^{m\times p}$
defined in (\ref{87}) can be rewritten as~:
\begin{equation}
\label{70}
T_{a}(x) \;=\; -K^{-1}SH(x)\  .
\end{equation}
In the following we show that we can find a real number $k^*\geq 1$ such 
that, if $k$ is  strictly larger than $k^*$, then the triple
$(kA,T_{a},B_{1m})$ satisfies all the assumptions of Theorem \ref{ExistenceObsAp}~:
\begin{remunerate}
\item
The forward completeness within $\OR$ 
is satisfied by assumption.
\item
(\ref{UCAp}) is satisfied since,
using (\ref{UnifInjH}) and the definition of $T_a$ in (\ref{70}),
we get, for each pair $(x_1,x_2 )$ in $\cl(\OR)^2$,
\begin{eqnarray*}
|x_1-x_2| & \leq &  \rho(|S^{-1} K| \, 
|K^{-1}SH(x_1)-K^{-1}SH(x_2)|)
\  ,
\\
& \leq & \rho(|S^{-1} K| \, |T_a(x_1)-T_a(x_2)|)
\ .
\end{eqnarray*}
\item
Let the function $\mathfrak E 
:\cl(\OR)\rightarrow\RR^{m\times p} $ be defined as~:
$$
\mathfrak{E}(x) \;=\; -(kA)^{-m}B_{1m}L_{f}^{m}b(h(x))\, .
$$
We have to show that this function satisfies
(\ref{Eq21}) and (\ref{EDFAp}).
Using (\ref{54}), we get, for each $x$ in $\OR$,
$$
\mathfrak{E}(x)   +B_{1m}b(h(x))= -K^{-1}SL_{f}H(x) + kAK^{-1}SH(x) 
= L_{f}T_a(x) -  kAT_a(x) 
\  .
$$
So, (\ref{Eq21}) does hold. Also, with (\ref{LipCond2}) and
$
k^{-m}  |K| \leq 1
$,
which holds for $k\geq 1$,
we get, for each $(x_1,x_2)$ in $\cl(\OR)^2$,
\begin{eqnarray*}
| \mathfrak{E}(x_1) \, - \, \mathfrak{E}(x_2) |
& =  & | (kA)^{-m}(B_{1m}(L_{f}^{m}(h(x_1)) \, - \, 
L_{f}^{m}(h(x_2)))) | 
\  ,
\\
& \leq & |(kA)^{-m}| \, |B_{1m}| \, L \, |H(x_1)-H(x_2)| 
\  ,
\\
& \leq & \frac{1}{\min_i |\lambda _i|^m} \, |B_{1m}| \, L \, |S^{-1}| \, 
|T_a(x_1)-T_a(x_2)| 
\  .
\end{eqnarray*}
Hence, (\ref{EDFAp}) is satisfied with
$
N\;=\; \frac{1}{\min_i |\lambda _i|^m} \, |B_{1m}| \, L \, |S^{-1}|
$
which does not depend on $k$.
\item
It remains to show that, by choosing $k$ sufficiently large, the 
constraint (\ref{RelLP}) is satisfied.
As $kA$ is a diagonal complex matrix, the inequality (\ref{RelLP}) is simply~:
$$
\frac{1}{\min_i |\lambda _i|^m} \, |B_{1m}| \, L \, |S^{-1}|\,  
 \frac{1}{k\,  (-\max_i\Re(\lambda _i))}\; <\; 1
\  .
$$
Clearly this inequality holds for all $k$ large enough.
\end{remunerate}
\subsection{Technical comments on section \ref{sec6}}
Due to space limitations, we give here only some hints on how the
results established for the case of completeness can be extended to
the case of boundedness observability.

The
introduction of $\gamma $ in the observer has mainly two consequences~:
\begin{remunerate}
\item
For the error convergence, $t$, in $\exp(At)$ in (\ref{85}) or
$\exp(-\varepsilon t)$ in (\ref{86}),
is replaced by the integral $\int_0^t \gamma (h(X(x,s)))ds$. If
$\sigma _{\RR^n}^+(x)=+\infty$, then,$\gamma $ being larger than $1$,
this integral goes to $+\infty $ as $t$ goes to $\sigma _{\RR^n}^+(x)$.
If, instead, $\sigma _{\RR^n}^+(x)$ is finite, then $\Vf(X(x,t))$ goes to
$+\infty $ as $t$ goes to $\sigma _{\RR^n}^+(x)$. From (\ref{ForUnb})
this is possible only if the above integral tends again to $+\infty $.
\item
The function $T$ given in (\ref{62}) is defined in terms
of the solutions $\bX(x,t)$ of the modified system (\ref{ModifSyst})
with $f_\gamma $ instead of $f$. So
we must show that this latter system shares
the backward $\OR $-distinguishability property of the original
system (\ref{3}).
This can be done by  associating,
to each $x$ in $\OR+\delta_d$, the function
$ \tau _x :
(\sigma _{\OR+\delta_d}^-(x),\sigma _{\OR+\delta_d}^+(x))
\rightarrow \RR $
defined as~:
$$
\tau _x (t)\;=\; \int_0^t\gamma(h(X(x,s)))ds
\  .
$$
It admits an inverse $\tau ^{-1}$ which is such that we have~:
$$
X(x,\tau^{-1}_x(t)) \;=\;  \bX(x,t)
\qquad \forall x \in \OR + \delta _d\; ,\ \forall t \in 
\tau_x(\sigma_{\OR+\delta_d}^-(x),\sigma_{\OR+\delta_d}^+(x))
\  .
$$
Then it is possible to prove that, if, for some pair $(x_1,x_2)$ in
$\OR^2$, we have~:
$$
h(\bX (x_1,t))\;=\; h(\bX (x_2,t))
\qquad \forall t\,  \in \,
(\bsigma _{\OR +\delta _d}^-(x_1),0] \cap(\bsigma _{\OR +\delta _d}^-(x_2),0]
\  ,
$$
then we have also~:
$$
\tau^{-1}_{x_1}(t)\;=\; \tau^{-1}_{x_2}(t)
\qquad \forall t\,  \in \,
(\bsigma _{\OR +\delta _d}^-(x_1),0] \cap(\bsigma _{\OR +\delta _d}^-(x_2),0]
\  .
$$
\end{remunerate}
\section{Conclusion}
We have stated sufficient conditions under which the extension to
non linear systems of the Luenberger observer, as it has been 
proposed 
by Kazantzis and Kravaris in \cite{Kazantzis-Kravaris}, can be used 
as long as the state to be 
observed remains in a given open set.
In doing so, we have exploited the fact, already mentioned in 
\cite{Astolfi-Praly,Krener-Xiao3},
that the observer 
proposed by Kreisselmeier and Engel in 
\cite{Kreisselmeier-Engel} is a possible way of implementing
the Kazantzis-Kravaris / Luenberger observer.

We have established that a sufficient (row) dimension of the dynamic 
system giving the observer is 2 + twice the dimension of the 
state to be observed. This is in agreement with many other results 
known on the generic number
of pieces of information to be extracted from
the output paths to be able to reconstruct the state.

We have also shown that it is sufficient to know only an 
approximation of a solution of a partial differential equation
which we need to solve to implement the observer. In this way, we 
have been able to make a connection with high gains observers.

{}finally, to get less restrictive sufficient conditions, we 
have 
found useful to modify the observer in a way which induces a time 
rescaling as already suggested in \cite{Astolfi-Praly}.

At this stage, our results are mainly of theoretical nature. 
They are concerned with existence. Several problems of prime 
importance for 
practice  remain to be addressed like type and speed of 
convergence. In these regards, the contribution of Rapaport and Maloum
in \cite{Rapaport-Maloum} is an important starting point.

Even for the purpose of showing the existence, we have to note that the 
conditions we have given can be strongly relaxed if an estimation
of the norm of the state is available. This idea has been exploited in
\cite{Astolfi-Praly} where a truly global observer has been proposed
under the assumption of global complete observability and 
unboundedness
observability.
\par\vspace{1em}
\Appendix
\section{Proof of Coron's Lemma \ref{Coron}}
The idea of the proof is to show that the set
$$
S\;=\; \bigcup_{x\in\Upsilon}
\left\{\Lambda =(\lambda _1,\ldots,\lambda _{n+1})\in \Omega^{n+1} 
 :\quad   g(x,\lambda_\ell) = 0
\quad \forall \ell\in\{1,\ldots,n+1\}\right\}
\ ,
$$
defined in (\ref{58}) is contained in a countable union of sets which have zero Lebesgue 
measure.

Given $(\underline x ,\underline \Lambda ,\epsilon)$ in
$\Upsilon\times \Omega ^{n+1} \times \RR_{+*}$, we denote by
$S_{\epsilon,\underline x, \underline \Lambda}$ the set~:
\begin{equation}
\label{DefS}
S_{\epsilon,\underline x, \underline \Lambda}
\;=\; \bigcup_{x\in \BR_\epsilon(\underline x)}
\left\{\Lambda \in 
\BR_\epsilon(\underline \Lambda)
\,   :\quad    g(x,\lambda_\ell)  =0\quad \forall \ell\in\{1,\ldots,n+1\}\right\}
\  .
\end{equation}
Assume for the time being that, for each pair
$(\underline x,\underline \Lambda)$ in $ \Upsilon\times \Omega ^{n+1}$, we
can find a strictly positive real number $\epsilon$ and a countable family of
$C^1$ functions 
$\sigma _{i}:\BR_{\epsilon}(\underline x)\rightarrow \Omega^{n+1}$, 
such that we have~:
\begin{equation}
\label{Objectif}
S_{\epsilon,\underline x, \underline \Lambda}
\subset \bigcup_{i\in \NN}\sigma _i(\BR_{\epsilon}(\underline x))
\  .
\end{equation}
The family
$\left(
\BR_{\epsilon}(\underline x)\times\BR_{\epsilon}(\underline \Lambda)
\right)_{(\underline x,\underline \Lambda) \in \Upsilon\times \Omega ^{n+1}
}$
is a covering of $\Upsilon\times \Omega ^{n+1}$ by open subsets.
From Lindel\"{o}f Theorem (see \cite[Lemma 4.1]{Boothby} for instance),
there exists a countable family
$\left\{(\underline x_j,\underline \Lambda_j)\right\}_{j\in\NN}$ such that
we have~:
$$
\Upsilon\times \Omega ^{n+1}\; \subset\; 
\bigcup_{j\in\NN}
\BR_{\epsilon_j}(\underline x_j)\times\BR_{\epsilon_j}(\underline 
\Lambda_j)
\  ,
$$
where $\epsilon_j$ denotes the $\epsilon$ associated to the pair
$(\underline x_j,\underline \Lambda_j)$.
With (\ref{Objectif}), it follows that we have~:
$$
S\; \subset\; \bigcup_{j\in\NN}
      \bigcup_{i\in \NN}\sigma _{i,j}(\BR_{\epsilon_j}(\underline x_j))
\  ,
$$
where $\sigma _{i,j}$ denotes the $i$th function $\sigma $ 
associated with the pair $(\underline x_j,\underline \Lambda_j)$.
The set $\sigma _{i,j}(\BR_{\epsilon_j}(\underline x_j))$
is the image, contained in $\CC^{n+1}$, a real manifold of dimension $2(n+1)$,
by a $C^1$ function of $\BR_{\epsilon_j}(\underline x_j)$, a real manifold of
dimension $2n$.
From a variation on Sard's Theorem (see \cite[Theorem 3, 
paragraphe 3]{Rham} for instance), this image
$\sigma _{i,j}(\BR_{\epsilon_j}(\underline x_j))$ has zero 
Lebesgue measure in $\CC^{n+1}$.
So $S$, being a countable union of such zero Lebesgue measure subsets,
has zero Lebesgue measure.

So all we have to do to establish Lemma \ref{Coron} is to prove the 
existence of $\varepsilon $ and the functions $\sigma _i$ satisfying (\ref{Objectif})
for each pair $(\underline x,\underline \Lambda) $in
$\Upsilon\times \Omega ^{n+1}$. For $\varepsilon $, we consider two
cases~:
\begin{remunerate}
\item
Consider a pair
$(\underline x,\underline \Lambda)$ such that
$g_j(\underline x, \underline \lambda_\ell)$ is non zero for some
component $\lambda _\ell$ of $\underline \Lambda$ and $g_j$ of $g$.
By continuity of $g_j$, we can
find a strictly positive real number $\epsilon $ such that
$g(x, \lambda_\ell)$ is also non zero for all $x$ in 
$\BR_\epsilon(\underline x)$ and
$\Lambda$ in $\BR_\epsilon(\underline \lambda)$.
In this case, the set
$S_{\epsilon,\underline x, \underline \Lambda}$ is empty. 
\item
Consider a pair $(\underline x,\underline \Lambda)$
such that $g(\underline x, \underline \lambda_\ell)$ is zero for 
each of the $n+1$ components of $\underline \lambda_\ell$ of
$\underline \Lambda$. From the assumption (\ref{Condg}),
for each $\ell$, we can find a component $g_{j_\ell}$ of $g$
and an integer $k_\ell$ satisfying~:
$$
\frac{\partial^i g_{j_\ell}}{\partial \lambda^i}(\underline x, \underline 
\lambda_\ell)
  \;=\; 0 \qquad \forall i \in\{0,\ldots,k_\ell-1\}
\qquad ,\qquad \qquad 
\frac{\partial^{k_\ell} g_{j_\ell}}{\partial \lambda^{k_\ell}}(\underline x, 
\underline \lambda_\ell)
   \; \neq\;  0
\  .
$$
In this case, following the Weierstrass Preparation Theorem
(see \cite[Theorem IV.1.1]{Golubitsky-Guillemin}\footnote{%
In \cite[Theorem IV.1.1]{Golubitsky-Guillemin}, this theorem is stated 
with 
the assumption that $g_j$ is holomorphic in both $x$ and $\lambda $. 
However, as far as $x$ is concerned, it can be seen in the proof of
this Theorem that we need only the implicit 
function theorem to apply. So continuous differentiability in $x$ for 
each $\lambda $ is enough.
} for instance), for each $\ell$ in $\{1,\ldots,n+1\}$, we know the existence of a 
strictly positive real number $\epsilon _\ell$,
a function
$q_\ell:\BR_{\epsilon_\ell}(\underline x)\times 
\BR_{\epsilon_\ell}(\underline \lambda _\ell)
\rightarrow\CC$, and $k_\ell$ $C^1$ functions
$a_j^\ell:\RR^{2n}\rightarrow\CC$ satisfying,
for all $(x,\lambda)$ in $\BR_{\epsilon_\ell}(\underline x) \times 
B_{\epsilon_\ell}(\underline\lambda_\ell)$,
\begin{equation}
\label{Weierstrass}
q_\ell(x,\lambda)\,  g_{j_\ell}(x,\lambda)
\;=\; (\lambda-\underline \lambda_\ell)^{k_\ell}
\;+\; \sum_{j=0}^{k_\ell-1}a_j^\ell(x)(\lambda-\underline 
\lambda_\ell)^j
\  .
\end{equation}
We choose the real number $\epsilon $, to be associated to
$(\underline x,\underline \Lambda)$ in the definition of $S_{\epsilon,\underline x, \underline
\Lambda}$, as~:
$$
\epsilon \;=\;  \inf_{\ell\in\{1,\ldots,n+1\}} \epsilon_\ell
\  .
$$
In the following
$P_\ell:\BR_\epsilon(\underline x)\times \CC\rightarrow\CC$
and $a^\ell (x):\BR_\epsilon(\underline x)\rightarrow \CC^{k_\ell }$
denote~:
$$
P_\ell(x,\lambda) 
\;=\;  (\lambda-\underline \lambda_\ell)^{k_\ell} 
\;+\; \sum_{j=0}^{k_\ell-1}a_j^\ell(x)(\lambda-\underline 
\lambda_\ell)^j
\quad ,\qquad  
a^\ell(x)\;=\; (a_0^\ell(x),\ldots,a_{k_\ell-1}^\ell(x))
\  .
$$
\end{remunerate}
\par\vspace{1em}
With this definition of $\varepsilon $, we have the following
implication, for $\Lambda $ in
$\BR_{\epsilon}(\underline \Lambda )$ and $x$ in $\BR_{\epsilon}(\underline
x)$,
\begin{equation}
\label{89}
\null\qquad \null
g(x,\lambda _\ell)= 0\quad \forall \ell\in\{1,\ldots,n+1\}
\quad \Rightarrow \quad 
(\lambda _\ell,a^\ell(x)) \in M^\ell\quad \forall \ell
\in\{1,\ldots,n+1\}
\end{equation}
where $M^\ell$ is the set~:
\begin{equation}
\label{Mi}
\null \qquad \null M^\ell \;=\; 
\left\{
 (\lambda,(b_0,\dots,b_{k_\ell-1}))\in \CC \times \CC^{k_\ell}
\: :\; (\lambda-\underline \lambda_\ell)^{k_\ell}
       + \sum_{j=0}^{k_\ell-1}b_j(\lambda-\underline \lambda_\ell)^j=0
\right\}
\end{equation}
Our interest in this set follows from
the following Lemma, which we prove later on,
\begin{lemma}
\label{StratLem}
Let $M$ be the set defined as~:
$$
M\;=\; 
\left\{
 (\lambda,b_0,\ldots,b_{k-1})\in \CC \times \CC^k
\: :\; \lambda^k
       + \sum_{j=0}^{k-1}b_j \lambda^j=0
\right\}\  .
$$
There exists a countable family $\{M_m\}_{m\in \NN}$ of regular submanifolds of 
$\CC^k$
and a countable family  of $C^1$ functions $\rho _m:M_m\rightarrow\CC$
such that we have the inclusion~:
\begin{equation}
\label{Strat2}
M\;\subset\;  \bigcup_{m \in \NN}
\bigcup_{b\in M_m}
\{(\rho _m(b),b)\}
\  .
\end{equation}
\end{lemma}

%
So, for each $\ell$ in $\{1,\ldots,n+1\}$
we have
a countable family $\{M_{m_\ell}^\ell\}_{m_\ell\in \NN}$ of regular submanifolds of 
$\CC^{k_\ell}$
and a countable family  of $C^1$
functions $\rho _{m_\ell}^\ell:M_{m_\ell}^\ell\rightarrow\CC$ such that,
for each $x$ in
$\BR_{\epsilon}(\underline x)$, if
$P_\ell(x,\lambda_\ell)$ is
zero, then there exists an integer
$m_\ell$ such that we have~:
\begin{equation}
\label{Pi}
a^\ell(x)\; \in\; M_{m_\ell}^\ell
\qquad ,\qquad \qquad 
\lambda _\ell\;=\; \rho _{m_\ell}^\ell (a^\ell(x))
\  .
\end{equation}
Hence, with (\ref{89}), if~:
$$
g(x,\lambda _\ell)\;=\; 0\qquad \forall \ell \in \{1,\ldots,n+1\}\  ,
$$
then there exists an $(n+1)$-tuple $
\mu =(m_1,\ldots,m_{n+1})$ of integers satisfying~:
$$
a^\ell(x)\; \in\; M_{m_\ell}^\ell
\  ,\quad  
\lambda _\ell\;=\; \rho _{m_\ell}^\ell (a^\ell(x))
\qquad \forall \ell \in \{1,\ldots,n+1\}\  .
$$
So, by letting~:
\begin{equation}
\label{60}
\null \qquad   \null 
S_{\epsilon,\underline x, \underline \Lambda}^\mu =
\bigcup_{
   \left\{x\in \BR_{\epsilon}(\underline x)
   \,   :\  
   a^\ell(x) \in M_{m_\ell}^\ell \  
\forall \ell \in \{1,\ldots,n+1\}\right\}
        }
\left\{
  (\rho _{m_1}^1(a^1(x)),\dots, \rho _{m_{n+1}} ^{n+1} (a^{n+1}(x))
\right\}
\end{equation}
we have established~:
\begin{equation}
\label{59}
S_{\epsilon,\underline x, \underline \Lambda}\; \subset\; 
\bigcup_{\mu \in \NN^{n+1}}
S_{\epsilon,\underline x, \underline \Lambda}^\mu 
\  .
\end{equation}
Comparing (\ref{Objectif}) with (\ref{59}) and using the definition 
(\ref{60}), we see that a candidate for the function $\sigma _i$ is~:
$$
\sigma _i(x) =
\left(
 \rho _{m_\ell}^\ell \left(R_{M_{m_\ell}^\ell}(a^\ell(x))\right)
\right)_{\ell\in\{1,\ldots,n+1\}}
$$
where $i$ happens to be the $(n+1)$-tuple $\mu $ and 
$R_{M_{m_\ell}^\ell}:\CC^{k_\ell}\rightarrow M_{m_\ell}^\ell$
is a ``restriction'' to $M_{m_\ell}^\ell$ since we 
have to consider only those $a^\ell(x)$ which are in $M_{m_\ell}^\ell $.
Finding such 
functions $R_{M_{m_\ell}^\ell}$ such that $\sigma _i$ is $C^1$ may 
not 
be possible. But, following \cite[Lemma 3.3]{Coron}, we know the 
existence, for each $\ell$, of a countable family of
$C^1$ functions
$R_{\nu }^\ell:\CC^{k_\ell}\rightarrow 
M_{m_\ell}^\ell$ such that
we get~:
$$
S_{\epsilon,\underline x, \underline \Lambda}^\mu
\; \subset \; \bigcup_{\nu \in \NN}
\left\{
\left(
 \rho _{m_\ell}^\ell\left( R_{\nu }^\ell(a^\ell(\BR_{\epsilon}(\underline x)))\right)
\right)_{\ell\in\{1,\ldots,n+1\}}
\right\}
\  .
$$
In other words the family of functions $\sigma _i$ is actually given 
by the family~:
$$
\sigma _{\mu ,\nu }\;=\; 
\left(
 \rho _{m_\ell}^\ell\circ R_\nu ^\ell \circ a_\ell
\right)_{\ell\in\{1,\ldots,n+1\}}
$$
i.e.~we have~:
\\[0.7em]\null \hfill $\displaystyle 
S_{\epsilon,\underline x, \underline \Lambda}\subset
\bigcup_{\mu \in \NN^{n+1}}\,  \bigcup_{\nu \in \NN}
\sigma _{\mu ,\nu }(\BR_{\epsilon}(\underline x))
\  .
$\null \ \null\hfill\Box\par\vspace{1em}
\noindent%
\textbf{Proof of Lemma \ref{StratLem} ~:}
The inclusion (\ref{Strat2}) says that we are looking for a 
covering of the set $M$ with some special structure. A covering 
easily found, but not having this special structure, is obtained by 
choosing a first
complex number $\lambda $, denoted $\lambda _1$, as well as $k-1$ 
other 
complex numbers $\lambda _j$. Then the $b_j$'s are given by the 
identity~:
\begin{equation}
\label{Rel}
\prod_{j=1}^{k} (\lambda -\lambda _j)\;=\; 
\lambda^k
       + \sum_{j=0}^{k-1}b_j \lambda^j
\  , \quad \lambda \in \CC\, .
\end{equation}
In other words if we denote by $\eta:\CC^k\rightarrow\CC^k$ the function 
which gives the $b_j$'s from the $\lambda_j$'s, we have~:
$$
\bigcup_{(\lambda_1,\dots,\lambda_k) \in \CC^{k}}
\{ ( \lambda_1,\eta(\lambda_1,\dots,\lambda_k))\} \subseteq M 
\  .
$$
Specifically, given the elementary 
symmetric functions $s_i$, sum of all the products of $i$
distinct $\lambda _j$'s,
$$
s_i\;=\; \lambda _1\ldots\lambda _{i-1}\lambda _i+
\lambda _1\ldots\lambda _{i-1}\lambda _{i+1}+\ldots
+\lambda _{k-i+1}\lambda _{k-i+2}\ldots\lambda _{k}\  ,
$$
the $b_j$'s are obtained as~:
$$
b_j\;=\; (-1)^{k-j}s_{k-j}\qquad \forall j\in\{1,\ldots,k-1\}
\  .
$$
Also the elementary symmetric functions are 
related to the sum of similar powers $\sigma _p$~:
$$
\sigma _p\;=\; \sum_{j=1}^k \lambda _j^p
\  ,
$$
via the Newton equations~:
$$
\sigma _i \,-\,  \sigma _{i-1}s_1 \,+\, \sigma _{i-2}s_2 \,+\, \dots 
\,+\, (-1)^{i-1}\sigma _1s_{k-1} \,+\, (-1) ^iis_i \;=\; 0
\  .
$$
The corresponding functions $(\sigma_p)\mapsto (s_i)$ and
$(s_i)\mapsto (b_j)$ are $C^\infty $ diffeomorphisms.

To obtain the result stated in the Lemma, we need to invert the 
function $\eta:(\lambda _\ell) \mapsto (b_j)$. 
This function is known to be an homeomorphism if the 
$\lambda _j$'s are defined up to permutations
(See \cite[Proposition 1.5.5]{Benedetti-Risler} for instance).
But unfortunately we cannot go beyond continuity of the inverse 
because of possible multiple roots. To round this problem, we choose 
the multiplicity $c_\ell$ of the root $\lambda _\ell$ so that the sum 
of the $c_\ell$'s is 
$k$. So, except if they are all $1$, some of them must be $0$.
Maybe after re-ordering, we can assume that each $c_1$ to $c_{q}$ is non zero and 
satisfy~:
$$
c_1\;+\; \ldots\;+\; c_{q}\;=\; k
\  .
$$
Then we choose $q$ different complex numbers $\varpi _\ell$ and we let~:
\\[0.7em]$\null \quad 
\lambda _1=\ldots =\lambda _{c_1}= \varpi _1
\  ,\  
\lambda _{c_0+1}=\ldots =\lambda _{c_1+c_2} = \varpi _2
\  ,\;  \ldots
$\hfill\null \\[0.3em]\null\hfill$
\ldots\  ,\; 
\lambda _{c_1+\ldots+c_{q-1}+1} =\ldots=\lambda _{c_1+\ldots+c_{q}}=\varpi _q
\  .
\quad \null $\\[0.7em]
This yields~:
\begin{equation}
\label{Stra1}
\sigma _p\;=\; \sum_{\ell=1}^{q} c_\ell \varpi _\ell^p
\  .
\end{equation}
We stress at this point that, to any $k$-tuple of $\lambda _\ell$'s in 
$\CC^k$, we can associate, maybe after a permutation $\theta $ of its 
components, such $q$-tuples of $c =(c_r)$ and
$\varpi =(\varpi _r)$, with $\varpi _i\neq 
\varpi _j$ if $i\neq j$.
It follows that the function $\eta$ can be decomposed as follows~:
$$
(\lambda_\ell)\; \underbrace{
\mapsto \; \theta (\lambda _\ell)\;
\mapsto \; (c,\varpi)\; 
\mapsto \; (\sigma _p)\;
\mapsto\; (s_i)\;
\mapsto \;}_{\eta}
(b_j)
\  .
$$
This way, given a permutation $\theta$ and a root multiplicity vector 
$c $, with no zero component, we have defined a 
function $\gamma:\CC^q\setminus\{\varpi _i = \varpi _j\}\,  \rightarrow\,  \CC^ k$
which maps the $\varpi _r$'s into the $b_j$'s~:
$$
\gamma : \varpi\; \mapsto \; (\sigma _p)\; \mapsto\; (s_i)\; 
\mapsto\; (b_j)
\  .
$$
This function has rank $q$. 
Indeed we know that the last two functions 
above are diffeomorphisms and, for the first one, we get from 
(\ref{Stra1})~:
$$
\frac{\partial \sigma_p}{\partial \varpi _r}\;=\; p\,   c_r \,  \varpi  _r^{p-1}
\  .
$$
Since the $p$'s and $c_r$'s are not zero, we see that the
matrix $\left(\frac{\partial \sigma_p}{\partial \varpi _r}\right)$ has full rank
since it has a Vandermonde 
structure and the $\varpi _r$'s are different.
Consequently, the jacobian matrix $\left(\frac{\partial b_j}{\partial \varpi_r}\right)$ of the function 
$\gamma$ has full rank $q\leq k$.
It follows from \cite[Theorem III,4.12, Theorem III, 5.5]{Boothby}
that, for each $q$-tuple $\varpi $ in
$\CC^q\setminus\{\varpi _i = \varpi _j\}$, there exists a strictly 
positive real number $\epsilon(\varpi )$ such that
\begin{itemize}
\item
$\BR _{\epsilon( \varpi)}( \varpi)$ is a subset of
$\CC^q\setminus\{\varpi _i = \varpi _j\}$,
\item
$\gamma \left(\BR _{\epsilon(\varpi))}( \varpi )\right)$
is a regular submanifold of (the real manifold) $\CC^k$,
\item
the restriction of $\gamma:
\BR _{\epsilon(\varpi)}(\varpi) \rightarrow
\gamma \left(\BR _{\epsilon(\varpi)}( \varpi )\right)$
is a diffeomorphism.
\end{itemize}
We denote by $\gamma ^{-1}$ the ``inverse'' function.

The family
$\left\{\BR _{\epsilon(\varpi)}( \varpi )\right\}_
    { \varpi \in \CC^q \setminus \left\{\varpi _i = \varpi _j\right\} }$ 
is a covering by open subsets of $\CC^{q}\setminus \{\varpi _i = \varpi _j\}$.
So there exists a countable family
$( \varpi^i)_{i\in\NN}$ such that the family
$\left\{\BR _{\epsilon( \varpi^i)}( \varpi^i)\right\}_{i\in \NN}$ is  covering
by open subsets of $\CC^{q} \setminus \{\varpi _i = \varpi _j\}$. Moreover,
since, to each $k$-tuple $(b_j)$ in $\CC^k$, we can associate a pair
$(c,\varpi)$ with $\varpi _i\neq \varpi _j$ if $i\neq j$,
any  such $k$-tuple $(b_j)$ is in at least one set
$\gamma \left(\BR _{\epsilon( \varpi^i)}( \varpi^i)\right)$.
So, since the number of
permutations $\theta$ in $\CC^k$ and multiplicity vectors 
$c $ is finite,
with varying $i$ and $q$, we get a countable family $\{M_m\}_{m\in\NN}$
of regular submanifold of $\CC^k$ defined as~:
$$
M_{m} \;:=\; \gamma \left(\BR _{\epsilon( \varpi^i)}( \varpi^i)\right)
$$
and a countable family of $C^1$ functions $\rho _m$ defined as~:
$$
\begin{array}{lllll}
&\gamma^{-1}& & (c ,\theta )&\\
\rho_m \; : \; (b_0,\dots,b_{k-1})\in M_m  & \mapsto &  
\varpi \in  \BR _{\epsilon( \varpi^i)}( \varpi^i)& \mapsto&
(\lambda_1,\dots,\lambda_k)\in \CC^k
\  .
\end{array}
$$
Each $b=(b_0,\dots,b_{k-1})$ in $\CC^k$ is in least one $M_m$ and we have~:
$$
\eta(\rho_m(b)) = b \, ,
$$
Our result follows then from~:
$$
\CC^k \; = \; \bigcup_{m \in \NN} \bigcup_{ b \in M_m}\{\rho_m(b)\}
\  .
$$

\textbf{Acknowledgements.}
This paper was initiated by many discussions with Dr. Alessandro 
Astolfi about the gap between the Kazantzis-Kravaris / Luenberger
observer, introduced as a local object in
\cite{Kazantzis-Kravaris}, and the global observer in
\cite{Astolfi-Praly}. We are also very grateful to Dr. Jean-Michel Coron who was our guide
all along our way in the quite involved technicalities needed in 
the proof of most of our results.

\end{document}